\documentclass{article}

\usepackage{amsmath,amssymb,amsfonts,amsthm,mdframed,graphicx}
\usepackage{color, colortbl}
\usepackage{xcolor}
\usepackage{array}
\usepackage{multirow}
\usepackage{enumitem}
\usepackage{url}

\newtheorem{theorem}{Theorem}
\newtheorem{lemma}{Lemma}

\newtheorem{proposition}[theorem]{Proposition}

\usepackage{algorithm}
\usepackage{algpseudocode}

\usepackage{tikz}
\usetikzlibrary{arrows.meta}
\tikzset{>={Latex[width=0.8mm,length=0.8mm]}}
\usetikzlibrary{positioning}
\usetikzlibrary{arrows}
\usepackage{tikz,amsmath, amssymb,bm,color}

\title{Algorithms for Short Human-Verifiable Proofs of Non-Word-Representability}

\author{Sergey Kitaev\footnote{Department of Mathematics and Statistics, University of Strathclyde, 26 Richmond Street, Glasgow G1, 1XH, United Kingdom. 
{\bf Email:} sergey.kitaev@strath.ac.uk.}\ \ and Junting Zhong\footnote{College of Science, Donghua University, Shanghai, 201620, P.R. China. {\bf Email:} jtzhong7@gmail.com}}

\begin{document}

\maketitle

\begin{abstract}
A graph is word-representable if it can be represented using alternation of letters in words. The class of word-representable graphs generalizes several important graph families and admits a characterization via semi-transitive orientations. Recognizing word-representability is an NP-complete problem, and a major bottleneck in the theory of word-representable graphs is providing efficient justifications of non-word-rep\-re\-sent\-ability of graphs, rather than merely referring to software verification.

Recently, Kitaev and Sun initiated the development of methods for the automatic generation of human-verifiable proofs of non-word-rep\-re\-sent\-ability. They introduced three algorithms and evaluated their performance on the 25 connected non-word-representable graphs on seven vertices, which constitute a fairly small sample set.

In this paper, we confirm their ranking of the algorithms by testing them on all non-word-representable graphs on eight and nine vertices. We further introduce three new algorithms, one of which significantly outperforms the previously best method on graphs with at most nine vertices, producing substantially shorter proofs. We apply this algorithm to shorten a known proof of non-word-representability of the Shrikhande graph and to produce a proof of non-word-representability of a certain subgraph of the Kneser graph $K(8,3)$. This subgraph is the first known example of a triangle-free non-word-representable graph, presented in the literature with reference only to a computational proof. As a by-product of our study, we correct an error in the repeatedly published list of the 25 connected non-word-representable graphs on seven vertices.  \\


\end{abstract}

\section{Introduction}

There is an extensive body of literature devoted to the theory of word-rep\-re\-sent\-able graphs (see, e.g.,~\cite{K17,KL15} and the references therein). The interest in this class stems from its connections to algebra, graph theory, computer science, combinatorics on words, and scheduling~\cite{KL15}. In particular, word-representable graphs generalize several fundamental graph classes, including \emph{circle graphs}, \emph{$3$-colorable graphs}, and \emph{comparability graphs}.

Two distinct letters $x$ and $y$ are said to \emph{alternate} in a word $w$ if after deleting from $w$ all letters except $x$ and $y$ the resulting word is of the form $xyxy\cdots$ or $yxyx\cdots$, of either even or odd length. A graph $G=(V,E)$ is \emph{word-representable} if and only if there exists a word $w$ over the alphabet $V$ such that, for distinct $x,y\in V$, the letters $x$ and $y$ alternate in $w$ if and only if $xy\in E$. In this case, we say that $w$ \emph{represents} $G$. A graph is word-representable if and only if each of its connected components is word-representable~\cite{KL15}. Therefore, without loss of generality, we may restrict our attention to connected graphs.

The smallest non-word-representable graph (in terms of the number of vertices) has six vertices; it is the wheel graph $W_5$ (obtained from Graph~1 in Figure~\ref{25-non-word-repr} by removing vertex~7). In contrast, there are $25$ connected non-word-representable graphs on seven vertices~\cite{KL15}. The numbers of such graphs on eight to eleven vertices are 929, 54957, 4880093, and 650856040, respectively.  This sequence took three years to compute and is listed as A290814 in the OEIS~\cite{oeis}.

An orientation of a graph is \emph{semi-transitive} if it is acyclic (i.e., there are no directed cycles) and satisfies the following condition: for every directed path
\[
v_0 \rightarrow v_1 \rightarrow \cdots \rightarrow v_k,
\]
either there is no edge between $v_0$ and $v_k$, or $v_i \rightarrow v_j$ is an edge for all $0 \le i < j \le k$.  An induced subgraph on vertices $\{v_0,v_1,\ldots,v_k\}$ of an oriented graph is a \emph{shortcut} if its orientation is acyclic and non-transitive, and it contains both the directed path $v_0 \rightarrow v_1 \rightarrow \cdots \rightarrow v_k$ and the edge $v_0 \rightarrow v_k$, which is referred to as the \emph{shortcutting edge}. Equivalently, a semi-transitive orientation can be defined as an acyclic orientation containing no shortcuts. A fundamental result in the area is the following characterization.

\begin{theorem}[{\cite{HKP16}}]\label{fundamental}
A graph is word-representable if and only if it admits a semi-transitive orientation.
\end{theorem}

The following known results are important for our paper.

\begingroup
\emergencystretch=2em
\begin{lemma}[\cite{KP}]\label{lemma} Suppose that an undirected graph $G$ has a cycle $C=x_1x_2\cdots x_mx_1$, where $m\geq 4$ and the vertices in $\{x_1,x_2,\ldots,x_m\}$ do not induce a clique in $G$.  If $G$ is oriented semi-transitively, and $m-2$ edges of $C$ are oriented in the same direction (i.e., from $x_i$ to $x_{i+1}$ or vice versa, where the index $m+1:=1$) then the remaining two edges of $C$ are oriented in the opposite direction.\end{lemma}
\endgroup
 
\begin{theorem}[\cite{KS}]\label{source-thm} Suppose that a graph $G$ is word-representable and that $v$ is a vertex of $G$. Then there exists a semi-transitive orientation of $G$ in which $v$ is a source (that is, all edges incident with $v$ are oriented away from $v$). Likewise, there exists a semi-transitive orientation of $G$ in which $v$ is a sink (that is, all edges incident with $v$ are oriented toward $v$). \end{theorem}

Recognizing word-representability of a graph is NP-complete~\cite{KL15}. A major practical difficulty in the theory of word-representable graphs is providing convincing and efficient proofs of non-word-representability, particularly when references to software (even publicly available and user-friendly tools such as~\cite{G}) are undesirable. In contrast, word-representability can be certified by exhibiting a semi-transitive orientation, which can be verified in polynomial time~\cite{KL15}.

\subsection{Proving non-word-representability} A comparability graph is an undirected graph that admits a transitive orientation. It is known (see~\cite{KP08}) that the neighbourhood of every vertex in a word-representable graph is a comparability graph, and recognizing comparability graphs can be done in polynomial time~\cite{G76}. Therefore, we obtain a polynomial-time test for non-word-representability of a graph~$G$: for each vertex, verify whether its neighbourhood induces a comparability graph; if a neighbourhood that is not a comparability graph is found, then $G$ is not word-representable. On the other hand, if all neighbourhoods are comparability graphs, the test is inconclusive, as the graph may or may not be word-representable~\cite{KL15}.

Thus, if the test above fails and no known non-word-representable induced subgraph is found in the given graph, then, prior to~\cite{KS}, the possible approaches to justify non-word-representability were the following:
\begin{itemize}[leftmargin=0.5cm]
\item[(a)] To produce some kind of ad-hoc (preferably elegant) argument that assumes the existence of a semi-transitive orientation and derives a contradiction. This approach is preferable, but such arguments are usually difficult to construct.

\item[(b)] To examine all $O(2^{\#\mbox{edges}})$ orientations and justify that none of them is semi-transitive. Symmetries can be exploited to reduce the search space; in particular, any given edge may be assumed to have a fixed orientation. A variation of this approach is implemented in some existing software~\cite{G,S}. It proceeds as follows: orient an edge $e_1$ in a given graph $G$, then select an as yet unoriented edge $e_2$ and branch on it, creating two copies of the partially oriented graph by orienting $e_2$ in the two possible ways; then branch on $e_3$, and so on. At each step, one ensures that no directed cycles or shortcuts are created (if they are, the corresponding branch is discarded). In any case, even if a computer were to output all these orientations (or the entire branching process) and indicate a directed cycle or a shortcut in each orientation, such a proof would be essentially unverifiable by a human, as it would typically be far too long.

\item[(c)] To examine all $O(\#\mbox{vertices}^2)$-length words in which each vertex label occurs the same number of times and to verify that none of them has the required alternation properties (if a graph with $n$ vertices is word-representable, then there exists a representing word of length at most $2n^2$~\cite{HKP16}). This approach requires checking a vast number of words. However, constraint programming can be used to speed up the process~\cite{Z}. In either case, it is impossible to convincingly demonstrate non-word-representability in this way without referring to software.

\item[(d)] Simply refer to publicly available software~\cite{G,S,Z}, which is less desirable, but still a valid option.
\end{itemize}

A game-changing approach that significantly increases the variety of non-word-representable graphs for which we can justify non-word-representability in a human-verifiable way was suggested in \cite{KS} by Kitaev and Sun. The basic idea in \cite{KS} is to modify the branching process using Lemma~\ref{lemma}, which helps avoid unnecessary branching through certain preprocessing steps. This modification enabled the development of software \cite{S} that can automatically generate human-verifiable proofs of non-word-representability.
The primary criterion for the efficiency of an algorithm is the number of “lines” (in the sense specified below) in the proof it produces; the secondary criterion is the running time of the algorithm. The latter is not discussed in this paper, as it is less relevant to our goals (the problem is NP-complete, and our focus is on proof length).
In fact, the approach in \cite{KS} is not novel: several papers, including \cite{KP}, use it or a simpler version of it (considering cycles of length 3 and 4 only) to justify non-word-representability. However, \cite{KS} is the first paper to discuss the automated search for human-verifiable proofs of non-word-representability of graphs, which led to the development of publicly available, user-friendly software \cite{S}. 

Theorem~\ref{source-thm} provides a further dramatic reduction in the size of proofs of non-word-representability and thereby further enlarges the class of non-word-representable graphs for which we can provide “short” proofs in a fully automatic manner. Indeed, selecting a vertex of maximum degree in a given non-word-representable graph, assuming it to be a source, and then running the branching process typically results in a much shorter proof. However, there are cases in which assuming a vertex of smaller degree to be a source leads to a shorter proof.

\subsection{Our contributions}

The paper \cite{KS} introduces three algorithms based on Lemma~\ref{lemma} to generate shorter proofs of non-word-representability. These algorithms were tested and ranked (by average performance) only on the 25 non-word-representable graphs on seven vertices, which is a rather small sample. We tested the algorithms on all non-word-representable graphs on eight and nine vertices and confirmed their ranking. Moreover, we propose three additional algorithms to address the problem, one of which significantly outperforms the best known algorithm so far, while the other two may be useful for particular graphs. In Section~\ref{algorithms}, we sketch all six algorithms and discuss their observed performance, and we provide pseudocode for the empirically best-performing one. We note that in our tests we chose not to use Theorem~\ref{source-thm}; employing this theorem and considering all possible assumptions about a source could lead to different results in the comparison of the algorithms.

\begin{center}
\begin{figure}
\begin{center}
\begin{tikzpicture}[scale=0.35]

\draw (0,0) node [scale=0.5, circle, draw](node5){5};
\draw (5,0) node [scale=0.5, circle, draw](node6){6};
\draw (6.54,4.75) node [scale=0.5, circle, draw](node2){2};
\draw (2.50,7.69) node [scale=0.5, circle, draw](node3){3};
\draw (-1.54,4.75) node [scale=0.5, circle, draw](node4){4};
\draw (2.5,1.5) node [scale=0.5, circle, draw](node7){7};
\draw(2.5,3.44)node [scale=0.5, circle, draw](node1){1};

\draw[->, >=stealth', line width=0.7pt] (node1)--(node2);
\draw[->, >=stealth', line width=0.7pt] (node1)--(node3);
\draw[->, >=stealth', line width=0.7pt] (node1)--(node4);
\draw[->, >=stealth', line width=0.7pt] (node1)--(node5);
\draw[->, >=stealth', line width=0.7pt] (node1)--(node6);
\draw[->, >=stealth', line width=0.7pt] (node1)--(node7);

\draw[->, >=stealth', line width=0.7pt] (node2)--(node3);
\draw[->, >=stealth', line width=0.7pt] (node2)--(node4);
\draw[->, >=stealth', line width=0.7pt] (node2)--(node6);
\draw[->, >=stealth', line width=0.7pt] (node2)--(node7);

\draw[->, >=stealth', line width=0.7pt] (node3)--(node4);
\draw[->, >=stealth', line width=0.7pt] (node3) to[out=-60,in=60] (node7);

\draw[->, >=stealth', line width=0.7pt] (node4)--(node7);

\draw[->, >=stealth', line width=0.7pt] (node5)--(node4);
\draw[->, >=stealth', line width=0.7pt] (node5)--(node6);
\draw[->, >=stealth', line width=0.7pt] (node5)--(node7);

\draw[->, >=stealth', line width=0.7pt] (node6)--(node7);


\end{tikzpicture}

\caption{The undirected version of this graph was assumed to be non-word-representable in the list of 25 non-word-representable graphs in \cite{KS}, although it is actually word-representable, as witnessed by the semi-transitive orientation shown in the figure. This graph is referred to as Graph~20 in \cite{KS}.}\label{a-graph}
\end{center}
\end{figure}
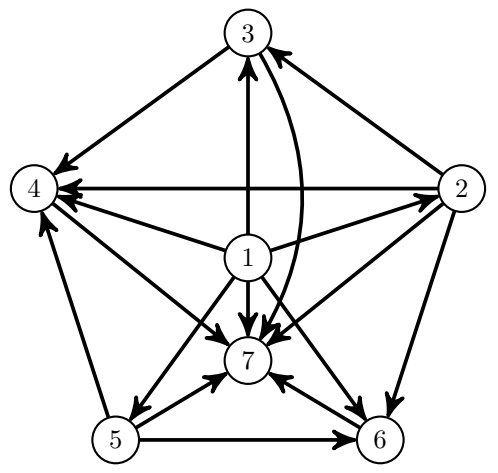
\end{center}


\begin{center}
\begin{figure}
\resizebox{\textwidth}{!}{%
\begin{tabular}{c c c}

\begin{tikzpicture}[scale=0.3]

\draw (0,12) node [scale=0.5, circle, draw](node1){1};
\draw (12,12) node [scale=0.5, circle, draw](node2){2};
\draw (12,0) node [scale=0.5, circle, draw](node3){3};
\draw (0,0) node [scale=0.5, circle, draw](node4){4};

\draw (6,11) node [scale=0.5, circle, draw](node5){5};
\draw (11,6) node [scale=0.5, circle, draw](node6){6};
\draw (6,1)node [scale=0.5, circle, draw](node7){7};
\draw (1,6)node [scale=0.5, circle, draw](node8){8};

\draw (6,9.5)node [scale=0.5, circle, draw](node9){9};
\draw (9.5,6)node [scale=0.38, circle, draw](node10){\large 10};
\draw (6,2.5)node [scale=0.38, circle, draw](node11){\large 11};
\draw (2.5,6)node [scale=0.38, circle, draw](node12){\large 12};

\draw (4.5,7.5)node [scale=0.38, circle, draw](node13){\large 13};
\draw (7.5,7.5)node [scale=0.38, circle, draw](node14){\large 14};
\draw (7.5,4.5)node [scale=0.38, circle, draw](node15){\large 15};
\draw (4.5,4.5)node [scale=0.38, circle, draw](node16){\large 16};

\draw (node1)--(node2)--(node3)--(node4)--(node1);

\draw (node1)--(node5);
\draw (node1)--(node9);
\draw (node1)--(node8);
\draw (node1)--(node12);

\draw (node2)--(node5);
\draw (node2)--(node9);
\draw (node2)--(node6);
\draw (node2)--(node10);

\draw (node3)--(node6);
\draw (node3)--(node7);
\draw (node3)--(node10);
\draw (node3)--(node11);

\draw (node4)--(node8);
\draw (node4)--(node7);
\draw (node4)--(node12);
\draw (node4)--(node11);

\draw (node8)--(node9)--(node6)--(node11)--(node8);
\draw (node5)--(node10)--(node7)--(node12)--(node5);

\draw (node13)--(node14)--(node15)--(node16)--(node13);

\draw (node5)--(node13);
\draw (node5)--(node14);
\draw (node11)--(node13);
\draw (node11)--(node14);

\draw (node9)--(node15);
\draw (node9)--(node16);
\draw (node7)--(node15);
\draw (node7)--(node16);

\draw (node8)--(node13);
\draw (node8)--(node16);
\draw (node10)--(node13);
\draw (node10)--(node16);

\draw (node6)--(node14);
\draw (node6)--(node15);
\draw (node12)--(node14);
\draw (node12)--(node15);

\node[text width=3cm] at (6.5,-1) 
    {Shrikhande graph};

\end{tikzpicture}

& 

\begin{tikzpicture}[scale=0.3]

\draw (0,0) node [scale=0.5, circle, draw](node4){4};

\draw (11,6) node [scale=0.5, circle, draw](node6){6};

\draw (1,6)node [scale=0.5, circle, draw](node8){8};

\draw (6,2.5)node [scale=0.38, circle, draw](node11){\large 11};
\draw (3,6)node [scale=0.38, circle, draw](node12){\large 12};

\draw (4.5,7.5)node [scale=0.38, circle, draw](node13){\large 13};
\draw (7.5,7.5)node [scale=0.38, circle, draw](node14){\large 14};
\draw (7.5,4.5)node [scale=0.38, circle, draw](node15){\large 15};
\draw (4.5,4.5)node [scale=0.38, circle, draw](node16){\large 16};

\draw (node4)--(node8);

\draw (node4)--(node12);
\draw (node4)--(node11);

\draw (node6)--(node11)--(node8);

\draw (node13)--(node14)--(node15)--(node16)--(node13);

\draw (node11)--(node13);
\draw (node11)--(node14);

\draw (node8)--(node13);
\draw (node8)--(node16);

\draw (node6)--(node14);
\draw (node6)--(node15);
\draw (node12)--(node14);
\draw (node12)--(node15);

\node[text width=3cm] at (7.5,-1) 
    {Subgraph $S$};

\end{tikzpicture}

&

\begin{tikzpicture}[scale=0.3]

\draw (0,0) node [scale=0.45, circle, draw](node4){2};

\draw (11,6) node [scale=0.45, circle, draw](node6){6};

\draw (1,6)node [scale=0.45, circle, draw](node8){3};

\draw (6,2.5)node [scale=0.45, circle, draw](node11){1};
\draw (3,6)node [scale=0.45, circle, draw](node12){7};

\draw (4.5,7.5)node [scale=0.45, circle, draw](node13){4};
\draw (7.5,7.5)node [scale=0.45, circle, draw](node14){5};
\draw (7.5,4.5)node [scale=0.45, circle, draw](node15){9};
\draw (4.5,4.5)node [scale=0.45, circle, draw](node16){8};

\draw (node4)--(node8);

\draw (node4)--(node12);
\draw (node4)--(node11);

\draw (node6)--(node11)--(node8);

\draw (node13)--(node14)--(node15)--(node16)--(node13);

\draw (node11)--(node13);
\draw (node11)--(node14);

\draw (node8)--(node13);
\draw (node8)--(node16);

\draw (node6)--(node14);
\draw (node6)--(node15);
\draw (node12)--(node14);
\draw (node12)--(node15);

\node[text width=3cm] at (7.5,-1) 
    {Subgraph $S$};

\end{tikzpicture}

\end{tabular}%
}
\caption{The Shrikhande graph and its minimal non-word-representable subgraph $S$ (with the original labelling and the labelling used in our proof)}\label{Shrikhande-pic}
\end{figure}
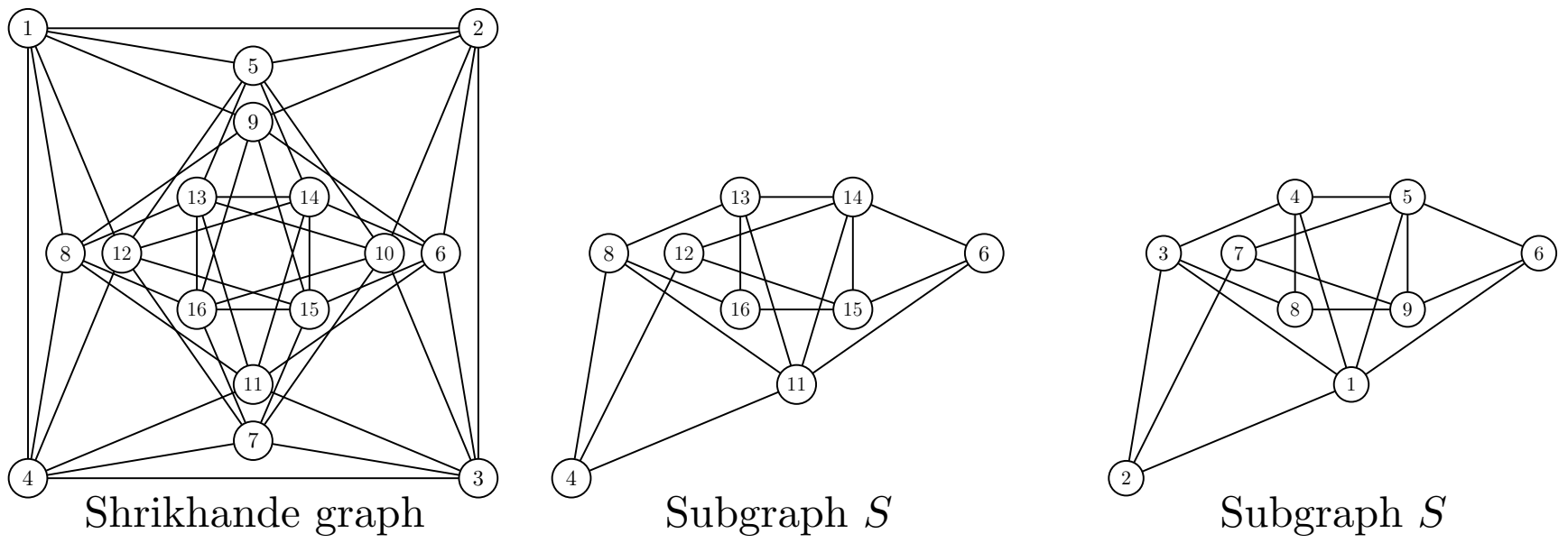
\end{center}

\begin{figure}
\begin{center}
\begin{tikzpicture}[node distance=1cm,auto,main node/.style={rectangle,rounded corners,draw,align=center}]

\node[main node] (1) {{\tiny {\bf 1}(357)}};
\node[main node] (2) [right of=1,xshift=6cm] {{\tiny {\bf 2}(468)}};

\node[main node] (3) [below of=1,yshift=-0.3cm] {{\tiny {\bf 3}(128)}};
\node[main node] (4) [below of=3,yshift=0.3cm] {{\tiny {\bf 4}(146)}};
\node[main node] (5) [below of=4,yshift=0.3cm] {{\tiny {\bf 5}(148)}};
\node[main node] (6) [below of=5,yshift=0.3cm] {{\tiny {\bf 6}(168)}};

\node[main node] (7) [below of=2,yshift=-0.3cm] {{\tiny {\bf 7}(217)}};
\node[main node] (8) [below of=7,yshift=0.3cm] {{\tiny {\bf 8}(235)}};
\node[main node] (9) [below of=8,yshift=0.3cm] {{\tiny {\bf 9}(237)}};
\node[main node] (10) [below of=9,yshift=0.3cm] {{\tiny {\bf 10}(257)}};

\node[main node] (11) [below of=6,yshift=-0.3cm] {{\tiny {\bf 11}(346)}};

\node[main node] (13) [below of=11,yshift=-0.3cm] {{\tiny {\bf 13}(258)}};
\node[main node] (14) [below of=13,yshift=0.3cm] {{\tiny {\bf 14}(458)}};

\node[main node] (12) [below of=10,yshift=-0.3cm] {{\tiny {\bf 12}(345)}};

\node[main node] (15) [below of=12,yshift=-0.3cm] {{\tiny {\bf 15}(167)}};
\node[main node] (16) [below of=15,yshift=0.3cm] {{\tiny {\bf 16}(367)}};

\path
(1) edge (2)
      edge (3)
      edge  [bend right=60] (4)
      edge [bend right=60] (5)
      edge [bend right=60] (6);
\path
(2) edge (7)
      edge  [bend left=60] (8)
      edge [bend left=60] (9)
      edge [bend left=60] (10);
\path
(3) edge [bend right=60] (11) edge (12) edge (16);
\path
(7) edge [bend left=60] (12) edge (11) edge (14);
\path
(4) edge (8) edge (9) edge (10) edge [bend right=60] (13);
\path
(5) edge (8) edge (9) edge (10) edge (16);
\path
(6) edge (8) edge (9) edge (10) edge (12);
\path
(8) edge [bend left=60] (15);
\path
(9) edge (14);
\path
(10) edge (11);
\path
(13) edge (15) edge (16) edge (11);
\path
(14) edge (15) edge (16);
\path
(12) edge (15);

\end{tikzpicture}
\caption{A minimal non-word-representable subgraph $S$ of the Kneser graph $K(8,3)$. Name of a vertex is in bold, and the corresponding 3-element subset is in parentheses.}\label{K83-fig}
\end{center}
\end{figure}
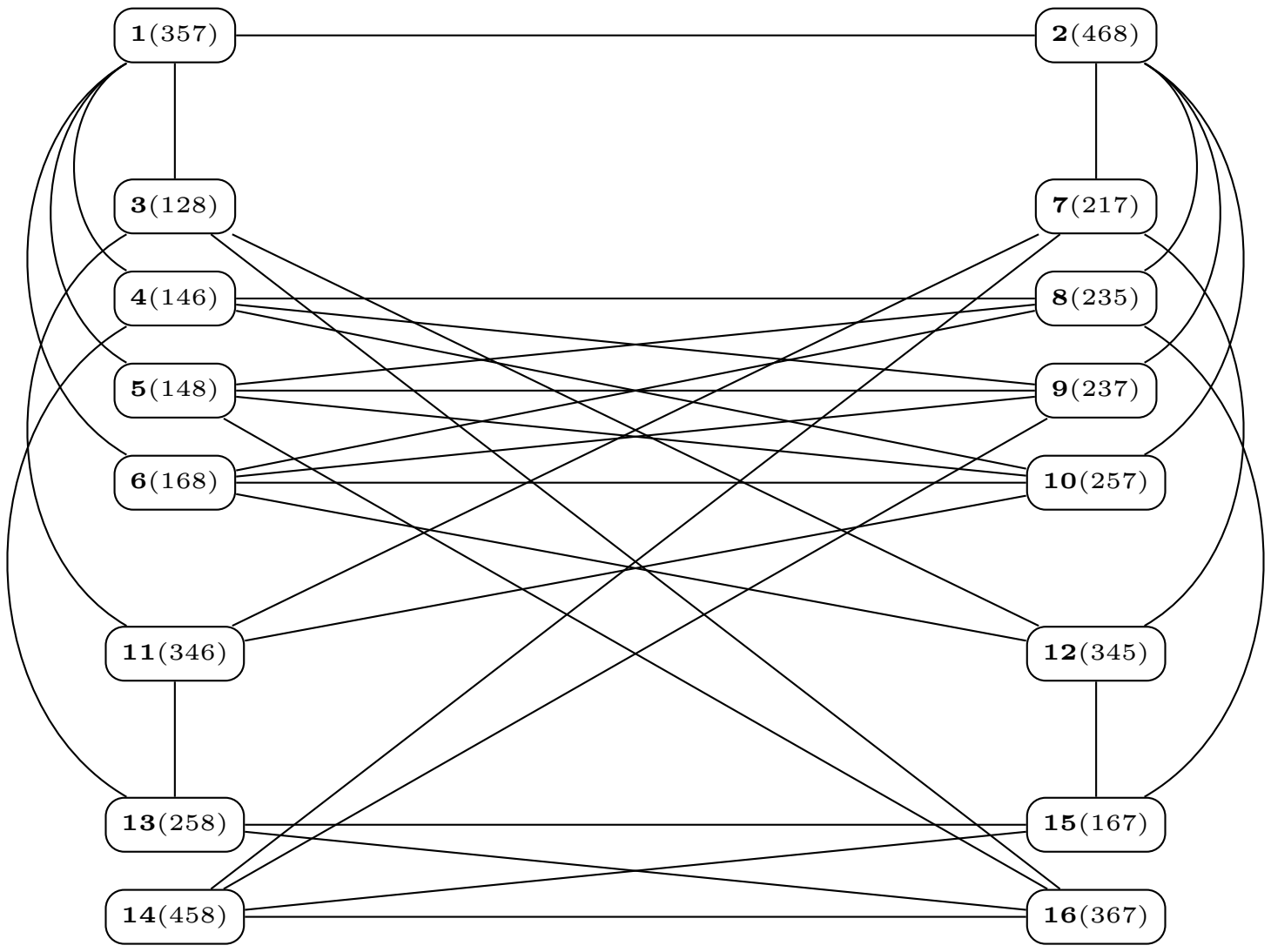

The paper \cite{KS} used software \cite{S} as a proof of concept to find automatically generated “short” proofs of non-word-representability for the \emph{Shrikhande graph} on 16 vertices and 48 edges (see Figure~\ref{Shrikhande-pic}), and the \emph{Clebsch graph} on 16 vertices and 40 edges. The latter graph has the even-cardinality subsets of $\{1,2,3,4,5\}$ as vertices, with two subsets adjacent if and only if their symmetric difference has cardinality four. While our best-performing algorithm (on average)  does not improve the length of the proof for the Clebsch graph, it reduces the proof length by one third for the Shrikhande graph. The publicly available software \cite{Zhong} can be used to confirm our results.

For any two integers $k \geq 1$ and $n \geq k$, the \emph{Kneser graph} $K(n,k)$ is the graph whose vertices correspond to the $k$-element subsets of the set $\{1,2,\ldots,n\}$, where two vertices are adjacent if and only if the corresponding sets are disjoint. Kneser graphs are vertex-transitive and edge-transitive, and they are named after Martin Kneser, who first investigated them in 1955. The subgraph $S$ of the Kneser graph $K(8,3)$ presented in \cite{KitSai} and displayed in Figure~\ref{K83-fig} is the first explicit example of a triangle-free non-word-representable graph; the existence of such graphs was established via Erd\H{o}s's theorem in \cite{HKP11}. A systematic study of the word-representability of triangle-free graphs was conducted in \cite{KP}. Unfortunately, the non-word-representability of $S$ was established in \cite{KitSai} only computationally (by referring to the software \cite{G}). We are now in a position to provide a human-verifiable proof of the non-word-representability of $S$; see Section~\ref{Shrikhande} for further details. 

\begin{figure}
\begin{center}
\resizebox{\textwidth}{!}{%
\begin{tabular}{c c c c c}

\begin{tikzpicture}[scale=0.25]

\draw (0,0) node [scale=0.5, circle, draw](node5){5};
\draw (5,0) node [scale=0.5, circle, draw](node6){6};
\draw (6.54,4.75) node [scale=0.5, circle, draw](node2){2};
\draw (2.50,7.69) node [scale=0.5, circle, draw](node3){3};
\draw (-1.54,4.75) node [scale=0.5, circle, draw](node4){4};
\draw (5,8.5) node [scale=0.5, circle, draw](node7){7};
\draw(2.5,3.44)node [scale=0.5, circle, draw](node1){1};

\draw (node6)--(node2)--(node3)--(node4)--(node5)--(node6)  ;
\draw (node1)--(node2);
\draw (node1)--(node3);
\draw (node1)--(node4);
\draw (node1)--(node5);
\draw (node1)--(node6);
\draw (node3)--(node7);
\draw (2.5,-1.8) node {Graph 1};
\end{tikzpicture} 

&
 
\begin{tikzpicture}[scale=0.25]

\draw (0,0) node [scale=0.5, circle, draw](node5){5};
\draw (5,0) node [scale=0.5, circle, draw](node6){6};
\draw (6.54,4.75) node [scale=0.5, circle, draw](node2){2};
\draw (2.50,7.69) node [scale=0.5, circle, draw](node3){3};
\draw (-1.54,4.75) node [scale=0.5, circle, draw](node4){4};
\draw (2.5,9.5) node [scale=0.5, circle, draw](node7){7};
\draw(2.5,3.44)node [scale=0.5, circle, draw](node1){1};

\draw (node6)--(node2)--(node3)--(node4)--(node5)--(node6)  ;
\draw (node1)--(node2);
\draw (node1)--(node3);
\draw (node1)--(node4);
\draw (node1)--(node5);
\draw (node1)--(node6);
\draw (node4)--(node7);
\draw (node2)--(node7);
\draw (2.5,-1.8) node {Graph 2};
\end{tikzpicture}
& 
\begin{tikzpicture}[scale=0.25]

\draw (0,0) node [scale=0.5, circle, draw](node5){5};
\draw (5,0) node [scale=0.5, circle, draw](node6){6};
\draw (6.54,4.75) node [scale=0.5, circle, draw](node2){2};
\draw (2.50,7.69) node [scale=0.5, circle, draw](node3){3};
\draw (-1.54,4.75) node [scale=0.5, circle, draw](node4){4};
\draw (2.5,9.5) node [scale=0.5, circle, draw](node7){7};
\draw(2.5,3.44)node [scale=0.5, circle, draw](node1){1};

\draw (node6)--(node2)--(node3)--(node4)--(node5)--(node6)  ;
\draw (node1)--(node2);
\draw (node1)--(node3);
\draw (node1)--(node4);
\draw (node1)--(node5);
\draw (node1)--(node6);
\draw (node4)--(node7);
\draw (node2)--(node7);
\draw (node3)--(node7);
\draw (2.5,-1.8) node {Graph 3};
\end{tikzpicture} 					
& 
\begin{tikzpicture}[scale=0.25]

\draw (0,0) node [scale=0.5, circle, draw](node5){5};
\draw (5,0) node [scale=0.5, circle, draw](node6){6};
\draw (6.54,4.75) node [scale=0.5, circle, draw](node2){2};
\draw (2.50,7.69) node [scale=0.5, circle, draw](node3){3};
\draw (-1.54,4.75) node [scale=0.5, circle, draw](node4){4};
\draw (2.5,9.5) node [scale=0.5, circle, draw](node7){7};
\draw(2.5,3.44)node [scale=0.5, circle, draw](node1){1};

\draw (node6)--(node2)--(node3)--(node4)--(node5)--(node6)  ;
\draw (node1)--(node2);
\draw (node1)--(node3);
\draw (node1)--(node4);
\draw (node1)--(node5);
\draw (node1)--(node6);
\draw (node4)--(node7);
\draw (node2)--(node7);
\draw    (node1) to[out=60,in=-60] (node7);
\draw (2.5,-1.8) node {Graph 4};

\end{tikzpicture}

&
 
\begin{tikzpicture}[scale=0.25]

\draw (0,0) node [scale=0.5, circle, draw](node5){5};
\draw (5,0) node [scale=0.5, circle, draw](node6){6};
\draw (6.54,4.75) node [scale=0.5, circle, draw](node2){2};
\draw (2.50,7.69) node [scale=0.5, circle, draw](node3){3};
\draw (-1.54,4.75) node [scale=0.5, circle, draw](node4){4};
\draw (6,8) node [scale=0.5, circle, draw](node7){7};
\draw(2.5,3.44)node [scale=0.5, circle, draw](node1){1};

\draw (node6)--(node2)--(node3)--(node4)--(node5)--(node6)  ;
\draw (node1)--(node2);
\draw (node1)--(node3);
\draw (node1)--(node4);
\draw (node1)--(node5);
\draw (node1)--(node6);
\draw (node1)--(node7);
\draw (2.5,-1.8) node {Graph 5};

\end{tikzpicture}
\\

\begin{tikzpicture}[scale=0.25]

\draw (0,0) node [scale=0.5, circle, draw](node5){5};
\draw (5,0) node [scale=0.5, circle, draw](node6){6};
\draw (6.54,4.75) node [scale=0.5, circle, draw](node2){2};
\draw (2.50,7.69) node [scale=0.5, circle, draw](node3){3};
\draw (-1.54,4.75) node [scale=0.5, circle, draw](node4){4};
\draw (6,8) node [scale=0.5, circle, draw](node7){7};
\draw(2.5,3.44)node [scale=0.5, circle, draw](node1){1};

\draw (node6)--(node2)--(node3)--(node4)--(node5)--(node6)  ;
\draw (node1)--(node2);
\draw (node1)--(node3);
\draw (node1)--(node4);
\draw (node1)--(node5);
\draw (node1)--(node6);
\draw (node1)--(node7);
\draw (node3)--(node7);

\draw (2.5,-1.8) node {Graph 6};

\end{tikzpicture}
& 
\begin{tikzpicture}[scale=0.25]

\draw (0,0) node [scale=0.5, circle, draw](node5){5};
\draw (5,0) node [scale=0.5, circle, draw](node6){6};
\draw (6.54,4.75) node [scale=0.5, circle, draw](node2){2};
\draw (2.50,7.69) node [scale=0.5, circle, draw](node3){3};
\draw (-1.54,4.75) node [scale=0.5, circle, draw](node4){4};
\draw (6,8) node [scale=0.5, circle, draw](node7){7};
\draw(2.5,3.44)node [scale=0.5, circle, draw](node1){1};

\draw (node6)--(node2)--(node3)--(node4)--(node5)--(node6)  ;
\draw (node1)--(node2);
\draw (node1)--(node3);
\draw (node1)--(node4);
\draw (node1)--(node5);
\draw (node1)--(node6);
\draw (node2)--(node7);
\draw (node3)--(node7);

\draw (2.5,-1.8) node {Graph 7};

\end{tikzpicture}
&  
\begin{tikzpicture}[scale=0.25]

\draw (0,0) node [scale=0.5, circle, draw](node5){5};
\draw (5,0) node [scale=0.5, circle, draw](node6){6};
\draw (6.54,4.75) node [scale=0.5, circle, draw](node2){2};
\draw (2.50,7.69) node [scale=0.5, circle, draw](node3){3};
\draw (-1.54,4.75) node [scale=0.5, circle, draw](node4){4};
\draw (6,8) node [scale=0.5, circle, draw](node7){7};
\draw(2.5,3.44)node [scale=0.5, circle, draw](node1){1};

\draw (node6)--(node2)--(node3)--(node4)--(node5)--(node6)  ;
\draw (node1)--(node2);
\draw (node1)--(node3);
\draw (node1)--(node4);
\draw (node1)--(node5);
\draw (node1)--(node6);
\draw (node2)--(node7);
\draw (node3)--(node7);
\draw (node1)--(node7);
\draw (2.5,-1.8) node {Graph 8};

\end{tikzpicture}
&

\begin{tikzpicture}[scale=0.25]

\draw (0,0) node [scale=0.5, circle, draw](node5){5};
\draw (5,0) node [scale=0.5, circle, draw](node6){6};
\draw (6.54,4.75) node [scale=0.5, circle, draw](node2){2};
\draw (2.50,7.69) node [scale=0.5, circle, draw](node3){3};
\draw (-1.54,4.75) node [scale=0.5, circle, draw](node4){4};
\draw (6,8) node [scale=0.5, circle, draw](node7){7};
\draw(2.5,3.44)node [scale=0.5, circle, draw](node1){1};

\draw (node6)--(node2)--(node3)--(node4)--(node5)--(node6)  ;
\draw (node1)--(node2);
\draw (node1)--(node3);
\draw (node1)--(node4);
\draw (node1)--(node5);
\draw (node1)--(node6);
\draw (node2)--(node7);
\draw (node3)--(node7);
\draw (node5)to[out=20,in=-100](node7);
\draw (2.5,-1.8) node {Graph 9};

\end{tikzpicture}
&
\begin{tikzpicture}[scale=0.25]

\draw (0,0) node [scale=0.5, circle, draw](node5){5};
\draw (5,0) node [scale=0.5, circle, draw](node6){6};
\draw (6.54,4.75) node [scale=0.5, circle, draw](node2){2};
\draw (2.50,7.69) node [scale=0.5, circle, draw](node3){3};
\draw (-1.54,4.75) node [scale=0.5, circle, draw](node4){4};
\draw (6,8) node [scale=0.5, circle, draw](node7){7};
\draw(2.5,3.44)node [scale=0.5, circle, draw](node1){1};

\draw (node6)--(node2)--(node3)--(node4)--(node5)--(node6)  ;
\draw (node1)--(node2);
\draw (node1)--(node3);
\draw (node1)--(node4);
\draw (node1)--(node5);
\draw (node1)--(node6);
\draw (node2)--(node7);
\draw (node3)--(node7);
\draw (node5)to[out=20,in=-100](node7);
\draw (node1)--(node7);
\draw (2.5,-1.8) node {Graph 10};

\end{tikzpicture}
\\

\begin{tikzpicture}[scale=0.25]

\draw (0,0) node [scale=0.5, circle, draw](node5){5};
\draw (5,0) node [scale=0.5, circle, draw](node6){6};
\draw (6.54,4.75) node [scale=0.5, circle, draw](node2){2};
\draw (2.50,7.69) node [scale=0.5, circle, draw](node3){3};
\draw (-1.54,4.75) node [scale=0.5, circle, draw](node4){4};
\draw (6,8) node [scale=0.5, circle, draw](node7){7};
\draw(2.5,3.44)node [scale=0.5, circle, draw](node1){1};

\draw (node6)--(node2)--(node3)--(node4)--(node5)--(node6)  ;
\draw (node1)--(node2);
\draw (node1)--(node3);
\draw (node1)--(node4);
\draw (node1)--(node5);
\draw (node1)--(node6);

\draw (node2)--(node7);
\draw (node3)--(node7);
\draw (node6)--(node7);
\draw (node1)--(node7);
\draw (2.5,-1.8) node {Graph 11};

\end{tikzpicture}
& 
\begin{tikzpicture}[scale=0.25]

\draw (0,0) node [scale=0.5, circle, draw](node5){5};
\draw (5,0) node [scale=0.5, circle, draw](node6){6};
\draw (6.54,4.75) node [scale=0.5, circle, draw](node2){2};
\draw (2.50,7.69) node [scale=0.5, circle, draw](node3){3};
\draw (-1.54,4.75) node [scale=0.5, circle, draw](node4){4};
\draw (6,8) node [scale=0.5, circle, draw](node7){7};
\draw(2.5,3.44)node [scale=0.5, circle, draw](node1){1};

\draw (node6)--(node2)--(node3)--(node4)--(node5)--(node6)  ;
\draw (node1)--(node2);
\draw (node1)--(node3);
\draw (node1)--(node4);
 
\draw (node1)--(node6);
\draw (node2)--(node7);
\draw (node3)--(node7);
\draw (node3)--(node5);
\draw (node3)--(node6);
\draw (node3)--(node1);
\draw (2.5,-1.8) node {Graph 12};

\end{tikzpicture}

&

\begin{tikzpicture}[scale=0.25]

\draw (0,0) node [scale=0.5, circle, draw](node5){5};
\draw (5,0) node [scale=0.5, circle, draw](node6){6};
\draw (6.54,4.75) node [scale=0.5, circle, draw](node2){2};
\draw (2.50,7.69) node [scale=0.5, circle, draw](node3){3};
\draw (-1.54,4.75) node [scale=0.5, circle, draw](node4){4};
\draw (2.5,1.5) node [scale=0.5, circle, draw](node7){7};
\draw(2.5,3.44)node [scale=0.5, circle, draw](node1){1};

\draw (node6)--(node2)--(node3)--(node4)--(node5)--(node6)  ;

\draw (node1)--(node3);

\draw (node1)--(node5);
\draw (node1)--(node6);
\draw (node7)--(node6);
\draw (node7)--(node2);
\draw (node7)--(node4);
\draw (node7)--(node5);
\draw (2.5,-1.8) node {Graph 13};

\end{tikzpicture}
& 
\begin{tikzpicture}[scale=0.25]

\draw (0,0) node [scale=0.5, circle, draw](node5){5};
\draw (5,0) node [scale=0.5, circle, draw](node6){6};
\draw (6.54,4.75) node [scale=0.5, circle, draw](node2){2};
\draw (2.50,7.69) node [scale=0.5, circle, draw](node3){3};
\draw (-1.54,4.75) node [scale=0.5, circle, draw](node4){4};
\draw (2.5,1.5) node [scale=0.5, circle, draw](node7){7};
\draw(2.5,3.44)node [scale=0.5, circle, draw](node1){1};

\draw (node6)--(node2)--(node3)--(node4)--(node5)--(node6)  ;

\draw (node1)--(node3);

\draw (node1)--(node5);
\draw (node1)--(node6);
\draw (node7)--(node6);
\draw (node7)--(node2);
\draw (node7)--(node4);
\draw (node7)--(node5);
\draw (node1)--(node4);
\draw (2.5,-1.8) node {Graph 14};

\end{tikzpicture}
& 
\begin{tikzpicture}[scale=0.25]

\draw (0,0) node [scale=0.5, circle, draw](node5){5};
\draw (5,0) node [scale=0.5, circle, draw](node6){6};
\draw (6.54,4.75) node [scale=0.5, circle, draw](node2){2};
\draw (2.50,7.69) node [scale=0.5, circle, draw](node3){3};
\draw (-1.54,4.75) node [scale=0.5, circle, draw](node4){4};
\draw (2.5,1.5) node [scale=0.5, circle, draw](node7){7};
\draw(2.5,3.44)node [scale=0.5, circle, draw](node1){1};

\draw (node6)--(node2)--(node3)--(node4)--(node5)--(node6)  ;

\draw (node1)--(node3);

\draw (node1)--(node5);
\draw (node1)--(node6);
\draw (node7)--(node6);
\draw (node7)--(node2);
\draw (node7)--(node4);
\draw (node7)--(node5);

\draw (node1)--(node4);
\draw (node2)--(node1);

\draw (2.5,-1.8) node {Graph 15};

\end{tikzpicture}
\\

\begin{tikzpicture}[scale=0.25]

\draw (0,0) node [scale=0.5, circle, draw](node5){5};
\draw (5,0) node [scale=0.5, circle, draw](node6){6};
\draw (6.54,4.75) node [scale=0.5, circle, draw](node2){2};
\draw (2.50,7.69) node [scale=0.5, circle, draw](node3){3};
\draw (-1.54,4.75) node [scale=0.5, circle, draw](node4){4};
\draw (2.5,1.5) node [scale=0.5, circle, draw](node7){7};
\draw(2.5,3.44)node [scale=0.5, circle, draw](node1){1};

\draw (node6)--(node2)--(node3)--(node4)--(node5)--(node6)  ;

\draw (node1)--(node3);

\draw (node1)--(node5);
\draw (node1)--(node6);
\draw (node7)--(node6);
\draw (node7)--(node2);
\draw (node7)--(node4);
\draw (node7)--(node5);

\draw (node1)--(node4);
\draw (node2)--(node1);
\draw (node7)--(node1);
\draw (2.5,-1.8) node {Graph 16};

\end{tikzpicture}
&

\begin{tikzpicture}[scale=0.25]

\draw (0,0) node [scale=0.5, circle, draw](node5){5};
\draw (5,0) node [scale=0.5, circle, draw](node6){6};
\draw (6.54,4.75) node [scale=0.5, circle, draw](node2){2};
\draw (2.50,7.69) node [scale=0.5, circle, draw](node3){3};
\draw (-1.54,4.75) node [scale=0.5, circle, draw](node4){4};
\draw (2.5,1.5) node [scale=0.5, circle, draw](node7){7};
\draw(2.5,3.44)node [scale=0.5, circle, draw](node1){1};

\draw (node6)--(node2)--(node3)--(node4)--(node5)--(node6)  ;

\draw (node1)--(node3);

\draw (node1)--(node5);
\draw (node1)--(node6);
\draw (node7)--(node6);
\draw (node7)--(node2);
\draw (node7)--(node4);
\draw (node7)--(node5);

\draw (node1)--(node4);
\draw (node2)--(node1);
\draw (node2)--(node4);
\draw (2.5,-1.8) node {Graph 17};

\end{tikzpicture}
&
\begin{tikzpicture}[scale=0.25]

\draw (0,0) node [scale=0.5, circle, draw](node5){5};
\draw (5,0) node [scale=0.5, circle, draw](node6){6};
\draw (6.54,4.75) node [scale=0.5, circle, draw](node2){2};
\draw (2.50,7.69) node [scale=0.5, circle, draw](node3){3};
\draw (-1.54,4.75) node [scale=0.5, circle, draw](node4){4};
\draw (2.5,1.5) node [scale=0.5, circle, draw](node7){7};
\draw(2.5,3.44)node [scale=0.5, circle, draw](node1){1};

\draw (node6)--(node2)--(node3)--(node4)--(node5)--(node6)  ;

\draw (node1)--(node3);

\draw (node1)--(node5);
\draw (node1)--(node6);
\draw (node7)--(node6);
\draw (node7)--(node2);
\draw (node7)--(node4);
\draw (node7)--(node5);

\draw (node1)--(node4);
\draw (node2)--(node1);
\draw (node7)to[out=60,in=-60](node3);
\draw (2.5,-1.8) node {Graph 18};

\end{tikzpicture}
&
\begin{tikzpicture}[scale=0.25]

\draw (0,0) node [scale=0.5, circle, draw](node5){5};
\draw (5,0) node [scale=0.5, circle, draw](node6){6};
\draw (6.54,4.75) node [scale=0.5, circle, draw](node2){2};
\draw (2.50,7.69) node [scale=0.5, circle, draw](node3){3};
\draw (-1.54,4.75) node [scale=0.5, circle, draw](node4){4};
\draw (2.5,1.5) node [scale=0.5, circle, draw](node7){7};
\draw(2.5,3.44)node [scale=0.5, circle, draw](node1){1};

\draw (node6)--(node2)--(node3)--(node4)--(node5)--(node6)  ;

\draw (node1)--(node3);

\draw (node1)--(node5);
\draw (node1)--(node6);
\draw (node7)--(node6);
\draw (node7)--(node2);
\draw (node7)--(node4);
\draw (node7)--(node5);
\draw (node7)--(node1);
\draw (node1)--(node4);
\draw (node2)--(node1);
\draw (node7)to[out=60,in=-60](node3);
\draw (2.5,-1.8) node {Graph 19};

\end{tikzpicture}
&
\begin{tikzpicture}[scale=0.25]

\draw(3.7,5)node [scale=0.5, circle, draw](node1){1};
\draw (2.5,1.5) node [scale=0.5, circle, draw](node2){2};
\draw (7.54,7.54) node [scale=0.5, circle, draw](node3){3};
\draw (7.54,0) node [scale=0.5, circle, draw](node4){4};
\draw (0,0) node [scale=0.5, circle, draw](node5){5};
\draw (0,7.54) node [scale=0.5, circle, draw](node6){6};
\draw (5.04,1.5) node [scale=0.5, circle, draw](node7){7};

\draw (node3)--(node4)--(node5)--(node6)--(node3);

\draw (node1)--(node2);
\draw (node1)--(node3);
\draw (node1)--(node4);
\draw (node1)--(node5);
\draw (node1)--(node6);
\draw (node1)--(node7);
\draw (node7)--(node2);
\draw (node7)--(node3);
\draw (node7)--(node4);
\draw (node2)--(node5);
\draw (node2)--(node6);

\draw (2.5,-1.8) node {\ \ \ \ \ Graph 20};

\end{tikzpicture}

\\

\begin{tikzpicture}[scale=0.25]

\draw (0,0) node [scale=0.5, circle, draw](node1){1};
\draw (7.54,0) node [scale=0.5, circle, draw](node2){2};
\draw (0,7.54) node [scale=0.5, circle, draw](node4){4};
\draw (7.54,7.54) node [scale=0.5, circle, draw](node3){3};
\draw (3,3) node [scale=0.5, circle, draw](node6){6};
\draw (3,1.2) node [scale=0.5, circle, draw](node5){5};
\draw(1.2,3)node [scale=0.5, circle, draw](node7){7};

\draw (node1)--(node2)--(node3)--(node4)--(node1)  ;
\draw (node4)--(node2);
\draw (node3)--(node6);
\draw (node1)--(node6);
\draw (node2)--(node6);
\draw (node7)--(node6);
\draw (node1)--(node7);
\draw (node1)--(node5);
\draw (node6)--(node5);
\draw (node2)--(node5);

\draw (4,-1.8) node {Graph 21};
\end{tikzpicture}

&

\begin{tikzpicture}[scale=0.25]

\draw (0,0) node [scale=0.5, circle, draw](node6){6};
\draw (4,6.84) node [scale=0.5, circle, draw](node5){5};
\draw (8,0) node [scale=0.5, circle, draw](node7){7};
\draw (4,2.33) node [scale=0.5, circle, draw](node1){1};
\draw (2,3.42) node [scale=0.5, circle, draw](node3){3};
\draw (6, 3.42) node [scale=0.5, circle, draw](node2){2};
\draw (4,0)node [scale=0.5, circle, draw](node4){4};

\draw (node6)--(node3)--(node5)--(node2)--(node7)--(node4)--(node6)  ;
\draw (node4)--(node2)--(node3)--(node4);
\draw (node1)--(node2);
\draw (node1)--(node3);
\draw (node1)--(node4);

\draw (4,-1.8) node {Graph 22};
\end{tikzpicture}

&

\begin{tikzpicture}[scale=0.25]

\draw (0,0) node [scale=0.5, circle, draw](node6){6};
\draw (4,6.84) node [scale=0.5, circle, draw](node5){5};
\draw (8,0) node [scale=0.5, circle, draw](node7){7};
\draw (4,2.33) node [scale=0.5, circle, draw](node1){1};
\draw (2,3.42) node [scale=0.5, circle, draw](node3){3};
\draw (6, 3.42) node [scale=0.5, circle, draw](node2){2};
\draw (4,0)node [scale=0.5, circle, draw](node4){4};

\draw (node6)--(node3)--(node5)--(node2)--(node7)--(node4)--(node6)  ;
\draw (node4)--(node2)--(node3)--(node4);
\draw (node1)--(node2);
\draw (node1)--(node3);
\draw (node1)--(node4);

\draw (node1)--(node5);
\draw (node1)--(node6);
\draw (node1)--(node7);
\draw (4,-1.8) node {Graph 23};
\end{tikzpicture} 
&
\begin{tikzpicture}[scale=0.25]

\draw (0,0) node [scale=0.5, circle, draw](node5){5};
\draw (5,0) node [scale=0.5, circle, draw](node6){6};
\draw (6.54,4.75) node [scale=0.5, circle, draw](node2){2};
\draw (2.50,7.69) node [scale=0.5, circle, draw](node3){3};
\draw (-1.54,4.75) node [scale=0.5, circle, draw](node4){4};
\draw (2.5,1.7) node [scale=0.5, circle, draw](node7){7};
\draw (2.5,3.44) node [scale=0.5, circle, draw](node1){1};


\draw (node1)--(node2);
\draw (node1)--(node3);
\draw (node1)--(node4);
\draw (node1)--(node5);
\draw (node1)--(node6);

\draw (node2)--(node4);
\draw (node7)--(node2);
\draw (node7)--(node4);
\draw (node7)--(node6);

\draw (node3)--(node2);
\draw (node3)--(node4);
\draw (node4)--(node5);
\draw (node5)--(node6);
\draw (node6)--(node2);

\draw (2.5,-1.8) node {Graph 24};

\end{tikzpicture}

&

\begin{tikzpicture}[scale=0.3]

\draw (0,0) node [scale=0.5, circle, draw](node6){6};
\draw (4,6.84) node [scale=0.5, circle, draw](node5){5};
\draw (8,0) node [scale=0.5, circle, draw](node7){7};
\draw (4,2.33) node [scale=0.5, circle, draw](node1){1};

\draw (2.6,2.84) node [scale=0.5, circle, draw](node3){3};
\draw (5.4, 2.84) node [scale=0.5, circle, draw](node2){2};
\draw (4,0.9)node [scale=0.5, circle, draw](node4){4};

\draw (node6)--(node5)--(node7)--(node6)  ;

\draw (node1)--(node2);
\draw (node1)--(node3);
\draw (node1)--(node4);

\draw (node1)--(node5);
\draw (node1)--(node6);
\draw (node1)--(node7);
\draw (node2)--(node7);
\draw (node3)--(node5);
\draw (node4)--(node6);
\draw (4,-1.8) node {Graph 25};
\end{tikzpicture} 
\end{tabular}%
}
\caption{All connected non-word-representable graphs on seven vertices}\label{25-non-word-repr}
\end{center}
\end{figure}
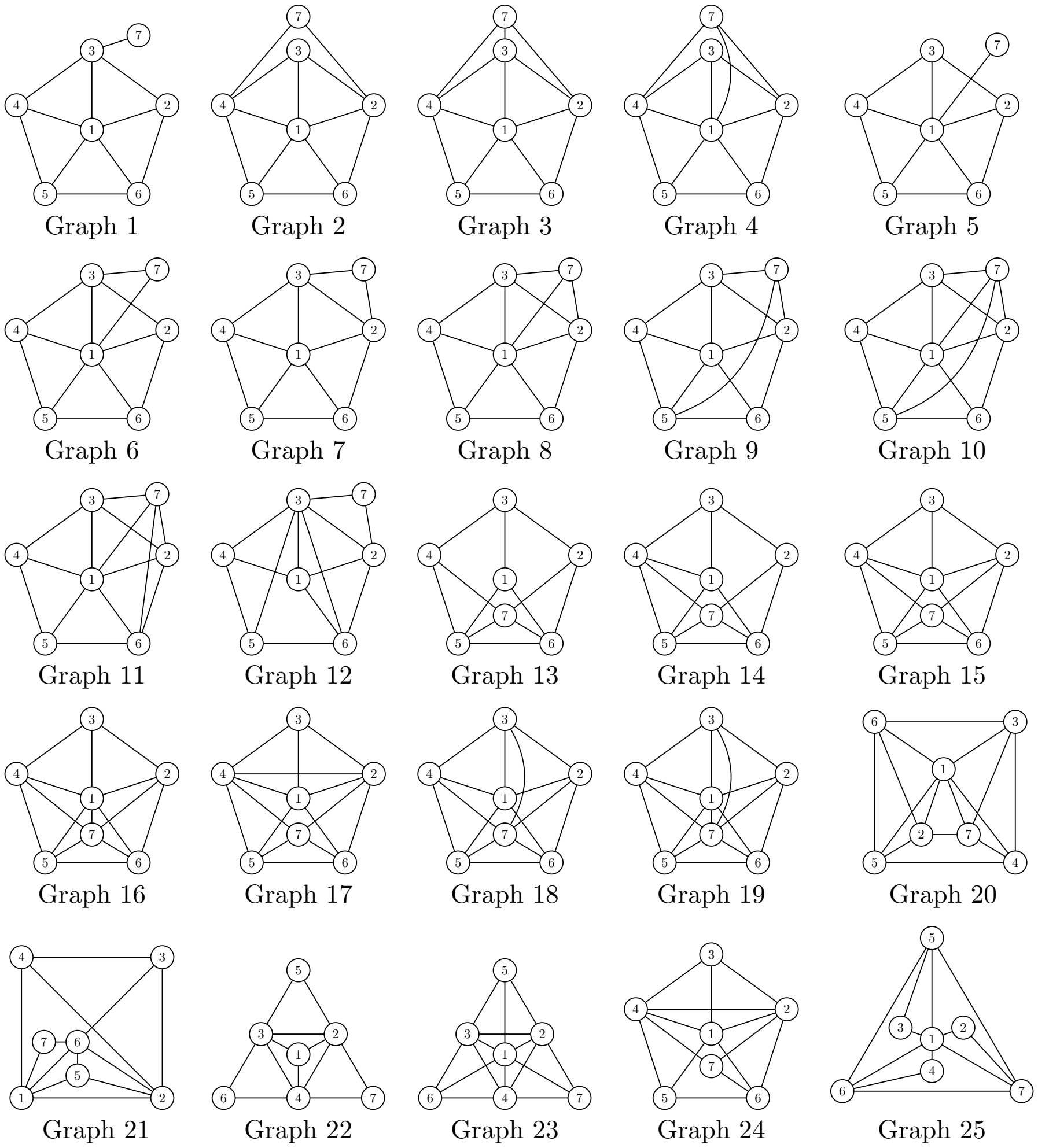

Finally, the list of all 25 connected non-word-representable graphs on seven vertices has been published multiple times, for example in \cite{KL15}. The paper \cite{KS} corrected two graphs in the list. However, we observed that the list presented in \cite{KS} still contains one incorrect graph. This graph, called Graph 20 in \cite{KS}, is the undirected version of the semi-transitively oriented graph shown in Figure~\ref{a-graph}. We leave it to the reader as a straightforward exercise to verify that the orientation in Figure~\ref{a-graph} is indeed semi-transitive. The correct list of the 25 connected non-word-representable graphs on seven vertices can be found in Figure~\ref{25-non-word-repr}, and a proof of non-word-representability of the updated Graph~20 follows from the easy-to-check fact that the neighborhood of vertex~1 is not a comparability graph. A proof for this graph, viewed as the line graph of the triangular book graph with three pages, is also given in~\cite[Lemma~10]{MK26}.

\section{Six algorithms for producing proofs}\label{algorithms}

In this section, we describe three known algorithms and three new algorithms for finding proofs of non-word-representability of graphs. All of these algorithms rely on the observation that the branching process should not involve any edges that do not belong to a cycle, since such edges can be oriented arbitrarily without affecting the semi-transitivity of the rest of the graph's orientation. Moreover, all three algorithms assume that, to produce a shorter proof, branching should be performed on edges that belong to many cycles, which tends to increase the number of applications of Lemma~\ref{lemma}. In any case, we begin by explaining the format of a proof of non-word-representability of a graph and give a small example below; see also~\cite{KS}. A completed exhaustive search is conclusive: either all branches close, proving non-word-representability, or a surviving terminal partial orientation can be completed semi-transitively. At such a terminal state, all edges belonging to non-clique cycles are already oriented. Orienting the remaining edges according to a topological ordering preserves acyclicity and cannot create a shortcut, since the path and shortcutting edge of any shortcut form a non-clique cycle. A search stopped before either outcome is established is inconclusive.

We write $uv$ for the undirected edge $\{u,v\}$ and $u\to v$ for an assigned direction. For lexicographic comparisons, vertex labels are ordered numerically and each edge $uv$ is represented by the pair $(u,v)$ with $u<v$, irrespective of any assigned direction. A cycle is represented by the lexicographically smallest vertex sequence among its cyclic shifts and their reversals.

\subsection{The format of a proof}

By a ``line'' of a proof, we mean a sequence of instructions that directs the orientation of a partially oriented graph and necessarily ends with a contradiction showing that this particular orientation branch will not produce a semi-transitive orientation. The idea is that if no branch produces a semi-transitive orientation, then the graph is not word-representable.

Each proof begins with assumptions on the orientations of certain edges, and there are four types of instructions:

\begin{itemize}[leftmargin=0.5cm]
\item ``B$a\rightarrow b$ (Copy $x$)'' means ``Branch on edge $ab$, orient the edge as $a\rightarrow b$, create a copy of the current version of the graph with the edge $ab$ oriented as $b\rightarrow a$, and call the new copy $x$; leave Copy $x$ aside and continue to follow the instructions.'' The instruction B occurs when no application of Lemma~\ref{lemma} is possible in the partially oriented graph.

\item ``MC$x$'' means ``Move to Copy $x$,'' where Copy $x$ is a partially oriented version of the graph created at some point in the branching process. This instruction is always followed by an oriented edge $a\rightarrow b$ to remind the reader of the directed edge obtained after branching.

\item ``O$a\rightarrow b$ (C$abc$)'' means orient the edge $ab$ as $a\rightarrow b$ in the cycle $abc$ to avoid a directed cycle. If the cycle is longer than a triangle, this represents an application of Lemma~\ref{lemma} to produce a cycle where all but two edges are oriented in one direction, and one of the remaining two edges is oriented in the opposite direction.

\item ``O$a\rightarrow b$ O$c\rightarrow d$ (C$xyz\cdots$)'' means that Lemma~\ref{lemma} is applied to cycle $xyz\cdots$ to create new directed edges $a\rightarrow b$ and $c\rightarrow d$.
\end{itemize}

\noindent
Each line ends with either ``S: $x$-$y$-$\cdots$-$z$,'' indicating a shortcut with directed path $x\to y\to\cdots\to z$ and shortcutting edge $x\to z$, or ``E: $x$-$y$-$\cdots$-$z$,'' indicating an orientation contradiction on the displayed cycle. Besides a directed cycle, this includes a non-clique cycle with directed path $x\to y\to\cdots\to z$ and the edge $xz$ still unoriented: either direction of $xz$ creates a directed cycle or a shortcut.

For example, consider the six-vertex wheel $W_5$ with hub~1 and rim $2$-$3$-$4$-$5$-$6$-$2$. By Theorem~\ref{source-thm}, we may assume that vertex~1 is a source, so $1\to i$ for $2\leq i\leq6$. The following is a two-line proof:

\begin{small}
\noindent
{\bf 1.} B2$\rightarrow$3 (Copy 2) O2$\rightarrow$6 (C1-3-2-6) O5$\rightarrow$6 (C1-2-6-5) O5$\rightarrow$4 (C1-4-5-6) O3$\rightarrow$4 (C1-3-4-5) S:1-2-3-4

\noindent
{\bf 2.} MC2 3$\rightarrow$2 O6$\rightarrow$2 (C1-3-2-6) O6$\rightarrow$5 (C1-2-6-5) O4$\rightarrow$5 (C1-4-5-6) O4$\rightarrow$3 (C1-3-4-5) S:1-4-3-2.
\end{small}

The first line assumes $2\to3$ and saves the alternative $3\to2$ as Copy~2. Each O instruction follows from Lemma~\ref{lemma} on the indicated four-cycle. The line ends with $1\to2\to3\to4$ and $1\to4$, a shortcut because $24$ is not an edge. The second line returns to Copy~2 and closes the other branch with the shortcut $1\to4\to3\to2$, again using the missing edge $24$. Thus both orientations of $23$ are excluded.

To print a proof for a single adjacency-matrix file \texttt{graph.txt}, run the software~\cite{Zhong} with
\begin{center}
\texttt{python3 proof\_lines.py --txt graph.txt --algorithm 6}.
\end{center}
The \texttt{--zip} option reports proof-length statistics for a collection of graphs instead.

\subsection{Algorithm 1 in \cite{KS}}

Algorithm 1 sorts the edges according to the number of cycles they belong to and then branches on an edge that belongs to the largest number of cycles (whenever branching is required). If there are two or more such edges, the branching choice is made lexicographically.

\subsection{Algorithm 2 (Best-Performing in \cite{KS})}

Algorithm 2 selects a cycle $C$ with the smallest number of unoriented edges. The unoriented edges in $C$ are then sorted, as in Algorithm 1, according to the number of cycles they belong to, and branching is performed on an edge that belongs to the largest number of cycles. If there are two or more such edges, the branching choice is made lexicographically.

\subsection{Algorithm 3 in \cite{KS}}

Algorithm 3 is similar to Algorithm 2, but it selects a cycle that has the largest number $N$ of edges oriented in the same direction. Among the cycles with the same value of $N$, Algorithm 3 selects a cycle $C$ with the smallest number of unoriented edges. Then, as in Algorithm 2, the unoriented edges in $C$ are sorted according to the number of cycles they belong to, and branching is performed on an edge that belongs to the largest number of cycles. If there are two or more such edges, the branching choice is made lexicographically.

\subsection{Algorithm 4 (New)}

Algorithm 4 attempts to improve the previously best-performing Algorithm~2 by keeping track of the length of the longest directed path in the current partially oriented graph (the idea being that longer directed paths may lead to more frequent applications of Lemma~\ref{lemma}).

Algorithm 4 selects a cycle $C$ with the smallest number of unoriented edges. Among the cycles with the same minimum number of unoriented edges, the choice is made lexicographically. If $C$ contains exactly one unoriented edge, branching is performed on this edge.

Otherwise, for each unoriented edge $e \in C$, Algorithm 4 temporarily orients $e$ in both possible directions and evaluates the resulting increase in the length of the longest directed path in the current partial orientation. For each edge, the maximum increase over the two possible orientations is recorded. The algorithm then selects an edge that maximizes this increase.

If several edges yield the same maximum increase, they are sorted in decreasing order according to the number of cycles they belong to, and branching is performed on the lexicographically smallest among them. If no edge produces a positive increase, Algorithm~4 falls back to the selection rule of Algorithm~3.

\subsection{Algorithm 5 (New)}

Algorithm 5 is another attempt to improve the previously best-performing Algorithm~2. Algorithm 5 selects a cycle $C$ with the smallest number of unoriented edges. Among the cycles with the same minimum number of unoriented edges, the choice is made lexicographically. The unoriented edges in $C$ are then sorted according to the following criteria:
\begin{enumerate}
    \item decreasing number of cycles the edge belongs to;
    \item decreasing sum of the degrees of its endpoints;
    \item lexicographic order.
\end{enumerate}

\noindent
Branching is performed on the first edge in this ordering.

\subsection{Algorithm 6 (New, Best-Performing)}

Algorithm 6 selects the branching edge using a deterministic one-step lookahead procedure. 
Unlike the previous algorithms, which rely only on local cycle information, 
Algorithm 6 explicitly simulates the immediate consequences of each candidate branching decision before selecting an edge.

The procedure consists of the following stages. \\[-3mm]

\noindent
{\bf 1. Selection of critical cycles.}
Among all non-clique cycles (as defined by Lemma~\ref{lemma}), we consider only those cycles that contain at least one unoriented edge. 
For each such cycle, we compute the number of its unoriented edges in the current partial orientation.

We then restrict our attention to the cycles with the \emph{minimum number of unoriented edges}. 
These cycles are the most constrained and are closest to being fully oriented. 
All further branching decisions are derived only from this subset of cycles.  \\[-3mm]

\noindent
{\bf 2. Construction of candidate edges.}
From the selected cycles, we collect all unoriented edges and remove duplicates. 
These edges form the candidate branching set. If this set contains only one edge, that edge is selected immediately.  \\[-3mm]

\noindent
{\bf 3. One-step simulation.}
For each candidate edge $e=uv$ with $u<v$, both possible orientations are examined:
\begin{itemize}
    \item Orient $u \to v$;
    \item Orient $v \to u$.
\end{itemize}

\noindent
For each orientation:
\begin{enumerate}
    \item The current state is copied and the copy is updated by orienting $e$ in the chosen direction.
    \item All forced propagation rules are applied: Lemma~\ref{lemma}, triangle rule, contradiction detection (directed cycles, shortcuts, or a non-clique cycle whose final unoriented edge would create one of these).
\end{enumerate}

If a contradiction is detected during propagation, the branch is considered immediately closed. 
Otherwise, we measure the remaining unresolved structure by counting the total number of unoriented cycle-edge incidences in the resulting state, with each edge counted once for each non-clique cycle containing it.

For each candidate edge, we therefore obtain:
\begin{itemize}
    \item The number of non-contradictory orientations,
    \item The worst remaining unresolved count among its valid orientations,
    \item The total remaining unresolved count over its valid orientations.
\end{itemize}

\noindent
{\bf 4. Edge comparison and selection.}
Candidate edges are compared using a fixed priority order:
\begin{enumerate}
    \item Prefer edges with fewer non-contradictory orientations (i.e., edges that close more branches immediately).
    \item Among those, prefer edges with smaller worst-case remaining unresolved structure.
    \item If still tied, prefer edges with smaller total remaining unresolved structure.
    \item If still tied, prefer edges that appear in more cycles.
    \item If still tied, prefer edges whose endpoints have larger degree sum.
    \item Finally, choose the lexicographically smallest edge, using the canonical pairs defined at the start of this section.
\end{enumerate}

The edge that is best under this ordered comparison is selected for branching.  \\[-3mm]

\noindent{\bf Interpretation.}
Algorithm 6 can be understood as a deterministic greedy lookahead rule. 
Instead of choosing an edge based solely on static structural properties, it evaluates the immediate impact of each branching decision and selects the edge with the best score in the stated priority order.

For the pseudocode, let $U_T(C)$ be the number of unoriented edges of cycle $C$ in state $T$, and let
\[
R(T)=\sum_{C\in\mathcal{C}} U_T(C).
\]
Thus $R(T)$ counts unoriented cycle-edge incidences, rather than distinct edges.

\setcounter{algorithm}{5}
\begin{algorithm}[p]
\caption{One-step lookahead branching}\label{best-alg-pseudocode}
\small
\begin{algorithmic}[1]
\Require $S$ \Comment{Current partial orientation state}
\Require $\mathcal{C}$ \Comment{Set of non-clique cycles}
\Require $\mathit{count}(e)$ \Comment{Number of cycles containing each edge}
\Ensure $e^*$ \Comment{Selected branching edge, or \texttt{None}}

\State $\mathit{critical\_cycles}\gets\{C\in\mathcal{C}\mid U_S(C)>0\}$
\If{$\mathit{critical\_cycles}=\varnothing$}
    \State \Return \texttt{None}
\EndIf

\State $t\gets\min\{U_S(C)\mid C\in\mathit{critical\_cycles}\}$
\State $\mathcal{C}^*\gets\{C\in\mathit{critical\_cycles}\mid U_S(C)=t\}$
\State $\mathit{candidates}\gets$ set of all unoriented edges in cycles in $\mathcal{C}^*$
\If{$|\mathit{candidates}|=1$}
    \State \Return the unique edge in $\mathit{candidates}$ as $e^*$
\EndIf

\For{each edge $e=uv$ in $\mathit{candidates}$, with $u<v$}
    \State $\mathit{unresolved}\gets0$
    \State $\mathit{worst\_residual}\gets0$
    \State $\mathit{residual\_sum}\gets0$

    \For{each direction $d\in\{u\to v,v\to u\}$}
        \State $\mathit{S\_copy}\gets\operatorname{clone}(S)$
        \State Update $\mathit{S\_copy}$ by orienting $e$ according to $d$
        \State Apply all forced propagation rules to $\mathit{S\_copy}$

        \If{a contradiction occurs}
            \State \textbf{continue} with the next direction in the inner loop
        \EndIf

        \State $\mathit{unresolved}\gets\mathit{unresolved}+1$
        \State $r\gets R(\mathit{S\_copy})$
        \State $\mathit{residual\_sum}\gets\mathit{residual\_sum}+r$
        \State $\mathit{worst\_residual}\gets\max(\mathit{worst\_residual},r)$
    \EndFor

    \State Assign $e$ the lexicographic score:
    \Statex \hspace{\algorithmicindent}$\begin{aligned}
    (&\mathit{unresolved},\mathit{worst\_residual},\mathit{residual\_sum},\\
     &-\mathit{count}(e),-\deg(u)-\deg(v),e).
    \end{aligned}$
\EndFor

\State $e^*\gets$ edge with minimal lexicographic score
\State \Return $e^*$
\end{algorithmic}
\end{algorithm}

The algorithm performs only one level of simulation and therefore remains computationally manageable while incorporating global structural information. Since this algorithm turns out to be the best-performing on average among all six algorithms, as shown by the experimental results in Section~\ref{ranking}, we present its pseudocode in Algorithm~\ref{best-alg-pseudocode}.

The following observation gives a formal guarantee for the lookahead rule.
\begin{proposition}\label{lookahead-guarantee}
Suppose that a candidate edge has $a$ trial orientations not closed by propagation and that each surviving trial state has at most $r$ unoriented cycle-edge incidences, taking $r=0$ when $a=0$. The search subtree following this branching choice has at most
\[
2-a+a2^r
\]
terminal branches. Algorithm~6 minimizes the number $a$, then the worst residual count $r$, and then the sum of the residual counts, among the candidate edges. Its branching rule preserves the soundness of the proof search.
\end{proposition}
\begin{proof}
The $2-a$ immediately closed orientations each contribute one terminal branch. In a surviving trial state, the number of distinct unoriented edges belonging to non-clique cycles is at most $r$. Each subsequent branching step orients one of these edges, so the resulting binary subtree has depth at most $r$ and at most $2^r$ terminal branches. This proves the bound. The minimization claim follows from the first three entries of the score. Finally, branching considers both directions, forced orientations are necessary in any semi-transitive extension, and branches are closed only by a directed cycle, a shortcut, or a non-clique cycle whose final unoriented edge would create one of these.
\end{proof}

This is a local guarantee; the residual count does not determine the eventual proof length. In particular, it does not establish that Algorithm~6 always gives shorter proofs than the other algorithms. The comparisons below are empirical.

\subsection{Ranking of algorithms}\label{ranking}

As noted in~\cite{KS}, making general statements about the efficiency of the algorithms, either in terms of the number of lines they produce or their time complexity, does not seem feasible. However, some indication of the algorithms' efficiency can be obtained by examining their performance on small non-word-representable graphs. Unfortunately, the first three algorithms were tested in~\cite{KS} only on non-word-representable graphs with up to seven vertices, which is a very small set. It was observed that Algorithms~2 and~3 have similar performance, with Algorithm~2 being slightly better (and hence used in the software~\cite{S}), while Algorithm~1 was much worse.

In Table~\ref{average-perf-table} we confirm the ranking of Algorithms~1--3 on graphs with up to nine vertices and extend the comparison by including Algorithms~4--6. On average, a clear winner is Algorithm~6, with Algorithms~2--5 performing similarly and Algorithm~1 being much worse (on nine vertices, the performance difference between the best and worst algorithms is a factor of about~2.09, and between the best and the next-best algorithms a factor of about~1.31).

\begin{table}[h]
\centering
\begin{tabular}{l|cccccc}
\hline
 & Alg.\ 1 & Alg.\ 2 & Alg.\ 3 & Alg.\ 4 & Alg.\ 5 & Alg.\ 6 \\
\hline
6 vertices & 10 & 8 & 8 & 8 & 8 & 8 \\
7 vertices & 12.880 & 8.080 & 8.280 & \cellcolor{gray!20!}  8.000 & 8.080 & 8.880 \\
8 vertices & 16.637 & 12.139 & 12.374 & 12.252 & 12.135 & \cellcolor{gray!20!}  9.688 \\
9 vertices & 22.350 & 14.046 & 14.218 & 14.256 & 14.045 & \cellcolor{gray!20!} 10.703 \\
\hline
\end{tabular}
\caption{Average proof length produced by the algorithms on graphs with six to nine vertices, without assuming a particular source vertex. The best-performing algorithm for graphs with seven or more vertices is highlighted in grey.}\label{average-perf-table}
\end{table}

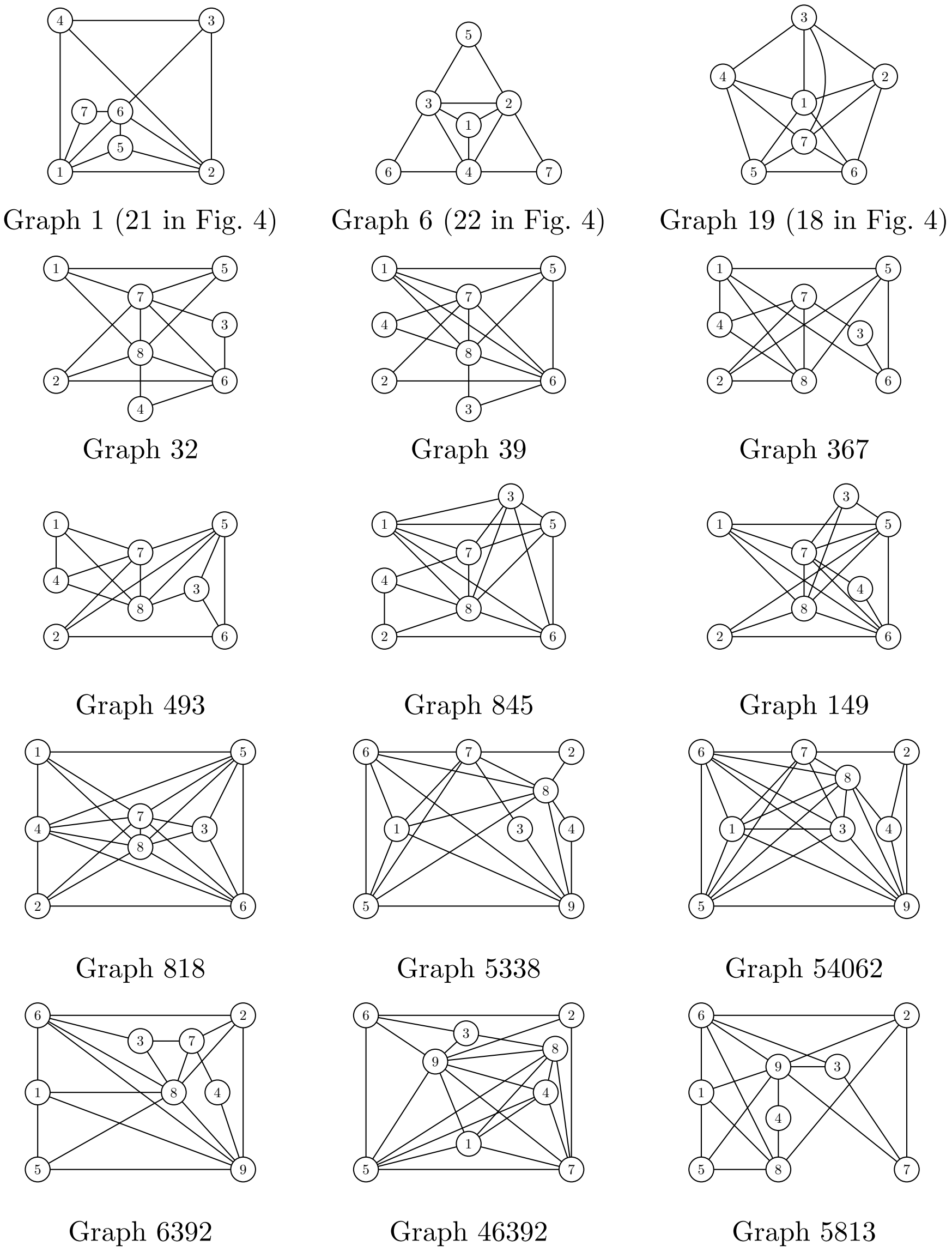
\begin{figure}[htbp]
\begin{center}
\resizebox{\textwidth}{!}{%
\begin{tabular}{c c c}
\begin{tikzpicture}[scale=0.25]

\draw (0,0) node [scale=0.5, circle, draw](node1){1};
\draw (7.54,0) node [scale=0.5, circle, draw](node2){2};
\draw (0,7.54) node [scale=0.5, circle, draw](node4){4};
\draw (7.54,7.54) node [scale=0.5, circle, draw](node3){3};
\draw (3,3) node [scale=0.5, circle, draw](node6){6};
\draw (3,1.2) node [scale=0.5, circle, draw](node5){5};
\draw(1.2,3)node [scale=0.5, circle, draw](node7){7};

\draw (node1)--(node2)--(node3)--(node4)--(node1)  ;
\draw (node4)--(node2);
\draw (node3)--(node6);
\draw (node1)--(node6);
\draw (node2)--(node6);
\draw (node7)--(node6);
\draw (node1)--(node7);
\draw (node1)--(node5);
\draw (node6)--(node5);
\draw (node2)--(node5);

\draw (4,-2.5) node {Graph 1 (21 in Fig.~\ref{25-non-word-repr})};
\end{tikzpicture}
& 
\begin{tikzpicture}[scale=0.25]

\draw (0,0) node [scale=0.5, circle, draw](node6){6};
\draw (4,6.84) node [scale=0.5, circle, draw](node5){5};
\draw (8,0) node [scale=0.5, circle, draw](node7){7};
\draw (4,2.33) node [scale=0.5, circle, draw](node1){1};
\draw (2,3.42) node [scale=0.5, circle, draw](node3){3};
\draw (6, 3.42) node [scale=0.5, circle, draw](node2){2};
\draw (4,0)node [scale=0.5, circle, draw](node4){4};

\draw (node6)--(node3)--(node5)--(node2)--(node7)--(node4)--(node6)  ;
\draw (node4)--(node2)--(node3)--(node4);
\draw (node1)--(node2);
\draw (node1)--(node3);
\draw (node1)--(node4);

\draw (4,-2.5) node {Graph 6 (22 in Fig.~\ref{25-non-word-repr})};
\end{tikzpicture}
& 
\begin{tikzpicture}[scale=0.25]

\draw (0,0) node [scale=0.5, circle, draw](node5){5};
\draw (5,0) node [scale=0.5, circle, draw](node6){6};
\draw (6.54,4.75) node [scale=0.5, circle, draw](node2){2};
\draw (2.50,7.69) node [scale=0.5, circle, draw](node3){3};
\draw (-1.54,4.75) node [scale=0.5, circle, draw](node4){4};
\draw (2.5,1.5) node [scale=0.5, circle, draw](node7){7};
\draw(2.5,3.44)node [scale=0.5, circle, draw](node1){1};

\draw (node6)--(node2)--(node3)--(node4)--(node5)--(node6)  ;

\draw (node1)--(node3);

\draw (node1)--(node5);
\draw (node1)--(node6);
\draw (node7)--(node6);
\draw (node7)--(node2);
\draw (node7)--(node4);
\draw (node7)--(node5);

\draw (node1)--(node4);
\draw (node2)--(node1);
\draw (node7)to[out=60,in=-60](node3);
\draw (2.5,-2.5) node {Graph 19 (18 in Fig.~\ref{25-non-word-repr})};
\end{tikzpicture}

\\

\begin{tikzpicture}[scale=0.35]

\draw (1,4)  node [scale=0.5, circle, draw](node1){1};
\draw (1,0)  node [scale=0.5, circle, draw](node2){2};
\draw (7,2)  node [scale=0.5, circle, draw](node3){3};
\draw (4,-1) node [scale=0.5, circle, draw](node4){4};
\draw (7,4)  node [scale=0.5, circle, draw](node5){5};
\draw (7,0)  node [scale=0.5, circle, draw](node6){6};
\draw (4,3)  node [scale=0.5, circle, draw](node7){7};
\draw (4,1)  node [scale=0.5, circle, draw](node8){8};

\draw (node1)--(node5);
\draw (node1)--(node7);
\draw (node1)--(node8);

\draw (node2)--(node6);
\draw (node2)--(node7);
\draw (node2)--(node8);

\draw (node3)--(node6);
\draw (node3)--(node7);

\draw (node4)--(node6);
\draw (node4)--(node8);

\draw (node5)--(node7);
\draw (node5)--(node8);

\draw (node6)--(node7);
\draw (node6)--(node8);

\draw (node7)--(node8);

\draw (4,-2.5) node {Graph 32};

\end{tikzpicture}
&
\begin{tikzpicture}[scale=0.35]

\draw (1,4)  node [scale=0.5, circle, draw](node1){1};
\draw (1,0)  node [scale=0.5, circle, draw](node2){2};
\draw (4,-1)  node [scale=0.5, circle, draw](node3){3};
\draw (1,2) node [scale=0.5, circle, draw](node4){4};
\draw (7,4)  node [scale=0.5, circle, draw](node5){5};
\draw (7,0)  node [scale=0.5, circle, draw](node6){6};
\draw (4,3)  node [scale=0.5, circle, draw](node7){7};
\draw (4,1)  node [scale=0.5, circle, draw](node8){8};

\draw (node1)--(node5);
\draw (node1)--(node6);
\draw (node1)--(node7);
\draw (node1)--(node8);

\draw (node2)--(node6);
\draw (node2)--(node7);

\draw (node3)--(node6);
\draw (node3)--(node8);

\draw (node4)--(node7);
\draw (node4)--(node8);

\draw (node5)--(node6);
\draw (node5)--(node7);
\draw (node5)--(node8);

\draw (node6)--(node7);
\draw (node6)--(node8);

\draw (node7)--(node8);

\draw (4,-2.5) node {Graph 39};

\end{tikzpicture}
&
\begin{tikzpicture}[scale=0.35]

\draw (1,4)  node [scale=0.5, circle, draw](node1){1};
\draw (1,0)  node [scale=0.5, circle, draw](node2){2};
\draw (6,1.7)  node [scale=0.5, circle, draw](node3){3};
\draw (1,2) node [scale=0.5, circle, draw](node4){4};
\draw (7,4)  node [scale=0.5, circle, draw](node5){5};
\draw (7,0)  node [scale=0.5, circle, draw](node6){6};
\draw (4,3)  node [scale=0.5, circle, draw](node7){7};
\draw (4,0)  node [scale=0.5, circle, draw](node8){8};

\draw (node1)--(node4);
\draw (node1)--(node5);
\draw (node1)--(node6);
\draw (node1)--(node8);

\draw (node2)--(node5);
\draw (node2)--(node7);
\draw (node2)--(node8);

\draw (node3)--(node6);
\draw (node3)--(node7);

\draw (node4)--(node7);
\draw (node4)--(node8);

\draw (node5)--(node6);
\draw (node5)--(node8);

\draw (node7)--(node8);

\draw (4,-2.5) node {Graph 367};

\end{tikzpicture}

\\

\begin{tikzpicture}[scale=0.35]

\draw (1,4)  node [scale=0.5, circle, draw](node1){1};
\draw (1,0)  node [scale=0.5, circle, draw](node2){2};
\draw (6,1.7)  node [scale=0.5, circle, draw](node3){3};
\draw (1,2) node [scale=0.5, circle, draw](node4){4};
\draw (7,4)  node [scale=0.5, circle, draw](node5){5};
\draw (7,0)  node [scale=0.5, circle, draw](node6){6};
\draw (4,3)  node [scale=0.5, circle, draw](node7){7};
\draw (4,1)  node [scale=0.5, circle, draw](node8){8};

\draw (node1)--(node4);
\draw (node1)--(node7);
\draw (node1)--(node8);

\draw (node2)--(node5);
\draw (node2)--(node6);
\draw (node2)--(node7);

\draw (node3)--(node5);
\draw (node3)--(node6);
\draw (node3)--(node8);

\draw (node4)--(node7);
\draw (node4)--(node8);

\draw (node5)--(node6);
\draw (node5)--(node7);
\draw (node5)--(node8);

\draw (node7)--(node8);

\draw (4,-2.5) node {Graph 493};

\end{tikzpicture}
&
\begin{tikzpicture}[scale=0.35]

\draw (1,4)  node [scale=0.5, circle, draw](node1){1};
\draw (1,0)  node [scale=0.5, circle, draw](node2){2};
\draw (5.5,5)  node [scale=0.5, circle, draw](node3){3};
\draw (1,2) node [scale=0.5, circle, draw](node4){4};
\draw (7,4)  node [scale=0.5, circle, draw](node5){5};
\draw (7,0)  node [scale=0.5, circle, draw](node6){6};
\draw (4,3)  node [scale=0.5, circle, draw](node7){7};
\draw (4,1)  node [scale=0.5, circle, draw](node8){8};

\draw (node1)--(node3);
\draw (node1)--(node5);
\draw (node1)--(node6);
\draw (node1)--(node7);
\draw (node1)--(node8);

\draw (node2)--(node4);
\draw (node2)--(node6);
\draw (node2)--(node8);

\draw (node3)--(node5);
\draw (node3)--(node6);
\draw (node3)--(node7);
\draw (node3)--(node8);

\draw (node4)--(node7);
\draw (node4)--(node8);

\draw (node5)--(node6);
\draw (node5)--(node7);
\draw (node5)--(node8);

\draw (node6)--(node8);

\draw (node7)--(node8);

\draw (4,-2.5) node {Graph 845};

\end{tikzpicture}
&
\begin{tikzpicture}[scale=0.35]

\draw (1,4)  node [scale=0.5, circle, draw](node1){1};
\draw (1,0)  node [scale=0.5, circle, draw](node2){2};
\draw (5.5,5)  node [scale=0.5, circle, draw](node3){3};
\draw (6,1.7) node [scale=0.5, circle, draw](node4){4};
\draw (7,4)  node [scale=0.5, circle, draw](node5){5};
\draw (7,0)  node [scale=0.5, circle, draw](node6){6};
\draw (4,3)  node [scale=0.5, circle, draw](node7){7};
\draw (4,1)  node [scale=0.5, circle, draw](node8){8};

\draw (node1)--(node5);
\draw (node1)--(node6);
\draw (node1)--(node7);
\draw (node1)--(node8);

\draw (node2)--(node5);
\draw (node2)--(node6);
\draw (node2)--(node8);

\draw (node3)--(node5);
\draw (node3)--(node7);
\draw (node3)--(node8);

\draw (node4)--(node6);
\draw (node4)--(node7);

\draw (node5)--(node6);
\draw (node5)--(node7);
\draw (node5)--(node8);

\draw (node6)--(node7);
\draw (node6)--(node8);

\draw (node7)--(node8);

\draw (4,-2.5) node {Graph 149};

\end{tikzpicture}

\\

\begin{tikzpicture}[scale=0.32]

\draw (0,6)  node [scale=0.5, circle, draw](node1){1};
\draw (0,0)  node [scale=0.5, circle, draw](node2){2};
\draw (6.5,3)  node [scale=0.5, circle, draw](node3){3};
\draw (0,3) node [scale=0.5, circle, draw](node4){4};
\draw (8,6)  node [scale=0.5, circle, draw](node5){5};
\draw (8,0)  node [scale=0.5, circle, draw](node6){6};
\draw (4,3.5)  node [scale=0.5, circle, draw](node7){7};
\draw (4,2.3)  node [scale=0.5, circle, draw](node8){8};

\draw (node1)--(node4);
\draw (node1)--(node5);
\draw (node1)--(node7);
\draw (node1)--(node8);

\draw (node2)--(node4);
\draw (node2)--(node6);
\draw (node2)--(node7);
\draw (node2)--(node8);

\draw (node3)--(node5);
\draw (node3)--(node6);
\draw (node3)--(node7);
\draw (node3)--(node8);

\draw (node4)--(node5);
\draw (node4)--(node6);
\draw (node4)--(node7);
\draw (node4)--(node8);

\draw (node5)--(node6);
\draw (node5)--(node7);
\draw (node5)--(node8);

\draw (node6)--(node7);
\draw (node6)--(node8);

\draw (node7)--(node8);

\draw (4,-2.5) node {Graph 818};

\end{tikzpicture}
&
\begin{tikzpicture}[scale=0.32]

\draw (1.2,3)  node [scale=0.5, circle, draw](node1){1};
\draw (8,6)  node [scale=0.5, circle, draw](node2){2};
\draw (6,3) node [scale=0.5, circle, draw](node3){3};
\draw (8,3) node [scale=0.5, circle, draw](node4){4};
\draw (0,0)  node [scale=0.5, circle, draw](node5){5};
\draw (0,6)  node [scale=0.5, circle, draw](node6){6};
\draw (4,6)  node [scale=0.5, circle, draw](node7){7};
\draw (7,4.5)  node [scale=0.5, circle, draw](node8){8};
\draw (8,0) node [scale=0.5, circle, draw](node9){9};

\draw (node1)--(node5);
\draw (node1)--(node6);
\draw (node1)--(node7);
\draw (node1)--(node8);
\draw (node1)--(node9);

\draw (node2)--(node7);
\draw (node2)--(node8);

\draw (node3)--(node7);
\draw (node3)--(node9);

\draw (node4)--(node8);
\draw (node4)--(node9);

\draw (node5)--(node6);
\draw (node5)--(node7);
\draw (node5)--(node8);
\draw (node5)--(node9);

\draw (node6)--(node7);
\draw (node6)--(node8);
\draw (node6)--(node9);

\draw (node7)--(node8);
\draw (node8)--(node9);

\draw (4,-2.5) node {Graph 5338};

\end{tikzpicture}
&
\begin{tikzpicture}[scale=0.32]

\draw (1.2,3)  node [scale=0.5, circle, draw](node1){1};
\draw (8,6)  node [scale=0.5, circle, draw](node2){2};
\draw (5.5,3) node [scale=0.5, circle, draw](node3){3};
\draw (7.3,3) node [scale=0.5, circle, draw](node4){4};
\draw (0,0)  node [scale=0.5, circle, draw](node5){5};
\draw (0,6)  node [scale=0.5, circle, draw](node6){6};
\draw (4,6)  node [scale=0.5, circle, draw](node7){7};
\draw (5.7,5)  node [scale=0.5, circle, draw](node8){8};
\draw (8,0) node [scale=0.5, circle, draw](node9){9};

\draw (node1)--(node3);
\draw (node1)--(node5);
\draw (node1)--(node6);
\draw (node1)--(node7);
\draw (node1)--(node8);
\draw (node1)--(node9);

\draw (node2)--(node4);
\draw (node2)--(node7);
\draw (node2)--(node9);

\draw (node3)--(node5);
\draw (node3)--(node6);
\draw (node3)--(node7);
\draw (node3)--(node8);
\draw (node3)--(node9);

\draw (node4)--(node8);
\draw (node4)--(node9);

\draw (node5)--(node6);
\draw (node5)--(node7);
\draw (node5)--(node8);
\draw (node5)--(node9);

\draw (node6)--(node7);
\draw (node6)--(node8);
\draw (node6)--(node9);

\draw (node7)--(node8);
\draw (node8)--(node9);

\draw (4,-2.5) node {Graph 54062};

\end{tikzpicture}

\\

\begin{tikzpicture}[scale=0.32]

\draw (0,3)  node [scale=0.5, circle, draw](node1){1};
\draw (8,6)  node [scale=0.5, circle, draw](node2){2};
\draw (4,5) node [scale=0.5, circle, draw](node3){3};
\draw (7,3) node [scale=0.5, circle, draw](node4){4};
\draw (0,0)  node [scale=0.5, circle, draw](node5){5};
\draw (0,6)  node [scale=0.5, circle, draw](node6){6};
\draw (6,5)  node [scale=0.5, circle, draw](node7){7};
\draw (5.3,3)  node [scale=0.5, circle, draw](node8){8};
\draw (8,0) node [scale=0.5, circle, draw](node9){9};

\draw (node1)--(node5);
\draw (node1)--(node6);
\draw (node1)--(node8);
\draw (node1)--(node9);

\draw (node2)--(node6);
\draw (node2)--(node7);
\draw (node2)--(node8);
\draw (node2)--(node9);

\draw (node3)--(node6);
\draw (node3)--(node7);
\draw (node3)--(node8);

\draw (node4)--(node7);
\draw (node4)--(node9);

\draw (node5)--(node8);
\draw (node5)--(node9);

\draw (node6)--(node8);
\draw (node6)--(node9);

\draw (node7)--(node8);
\draw (node8)--(node9);

\draw (4,-2.5) node {Graph 6392};

\end{tikzpicture}
&
\begin{tikzpicture}[scale=0.32]

\draw (4,1)  node [scale=0.5, circle, draw](node1){1};
\draw (8,6)  node [scale=0.5, circle, draw](node2){2};
\draw (3.9,5.3) node [scale=0.5, circle, draw](node3){3};
\draw (7,3) node [scale=0.5, circle, draw](node4){4};
\draw (0,0)  node [scale=0.5, circle, draw](node5){5};
\draw (0,6)  node [scale=0.5, circle, draw](node6){6};
\draw (8,0)  node [scale=0.5, circle, draw](node7){7};
\draw (7.37,4.7)  node [scale=0.5, circle, draw](node8){8};
\draw (2.7,4.2) node [scale=0.5, circle, draw](node9){9};

\draw (node1)--(node4);
\draw (node1)--(node5);
\draw (node1)--(node7);
\draw (node1)--(node8);
\draw (node1)--(node9);

\draw (node2)--(node6);
\draw (node2)--(node7);
\draw (node2)--(node9);

\draw (node3)--(node6);
\draw (node3)--(node8);
\draw (node3)--(node9);

\draw (node4)--(node5);
\draw (node4)--(node7);
\draw (node4)--(node8);
\draw (node4)--(node9);

\draw (node5)--(node6);
\draw (node5)--(node7);
\draw (node5)--(node8);
\draw (node5)--(node9);

\draw (node6)--(node9);

\draw (node7)--(node8);
\draw (node7)--(node9);
\draw (node8)--(node9);

\draw (4,-2.5) node {Graph 46392};

\end{tikzpicture}
&
\begin{tikzpicture}[scale=0.32]

\draw (0,3)  node [scale=0.5, circle, draw](node1){1};
\draw (8,6)  node [scale=0.5, circle, draw](node2){2};
\draw (5.3,4) node [scale=0.5, circle, draw](node3){3};
\draw (3,2) node [scale=0.5, circle, draw](node4){4};
\draw (0,0)  node [scale=0.5, circle, draw](node5){5};
\draw (0,6)  node [scale=0.5, circle, draw](node6){6};
\draw (8,0)  node [scale=0.5, circle, draw](node7){7};
\draw (3,0)  node [scale=0.5, circle, draw](node8){8};
\draw (3,4) node [scale=0.5, circle, draw](node9){9};

\draw (node1)--(node5);
\draw (node1)--(node6);
\draw (node1)--(node8);
\draw (node1)--(node9);

\draw (node2)--(node6);
\draw (node2)--(node7);
\draw (node2)--(node8);
\draw (node2)--(node9);

\draw (node3)--(node6);
\draw (node3)--(node7);
\draw (node3)--(node9);

\draw (node4)--(node8);
\draw (node4)--(node9);

\draw (node5)--(node8);
\draw (node5)--(node9);

\draw (node6)--(node8);
\draw (node6)--(node9);

\draw (node7)--(node9);

\draw (4,-2.5) node {Graph 5813};

\end{tikzpicture}

\end{tabular}%
}
\caption{Graphs in Table~\ref{worst-perf-table} that yield the longest proofs for the respective algorithms. The graph names come from the database~\cite{Zhong-graphs}. Note that the top three graphs also appear in Figure~\ref{25-non-word-repr}.}\label{15-non-word-repr}
\end{center}
\end{figure}

Interestingly, in the worst-case scenario (when the longest proofs are produced), Algorithm~6 significantly outperforms all other algorithms (on nine vertices, the performance difference ratio is between about~2.18 and~2.68), while Algorithm~1 is not the worst in the worst-case scenario on eight and nine vertices. We refer to Table~\ref{worst-perf-table} for the recorded observations and to Figure~\ref{15-non-word-repr}, where the graphs requiring long proofs are presented (their names come from the database~\cite{Zhong-graphs}).

\begin{table}[h]
\centering
\begin{tabular}{l|cl|cl|cl}
\hline
 &  Worst 7 & graphs &  Worst 8 & graphs &  Worst 9 & graphs \\
\hline
Alg.\ 1  &  28 & 1  &   56 & 32, 39, 367, 493 &   200 & 5338 \\
Alg.\ 2 &  12 & 1, 6 &  60 & 845 &   228 & 54062   \\
 Alg.\ 3 &  14 & 19 &  52 &  845 &    185 & 6392 \\
Alg.\ 4 &  12 & 1  &   52 &  845 &    188 &  46392  \\
Alg.\ 5 & 12 & 1, 6 & 60 & 845 &  228 & 54062 \\
Alg.\ 6 & 12 & 1 & 23 & 149, 818 &  85 & 5813 \\
\hline
\end{tabular}
\caption{Longest proofs produced by the implemented algorithms. ``Worst~X'' refers to the length of the longest proof output by an algorithm on a graph with $X$ vertices, and ``graphs'' refers to the graph(s) on which these proofs are achieved. The graphs can be found in Figure~\ref{15-non-word-repr} and are numbered as in the database~\cite{Zhong-graphs}.}
\label{worst-perf-table}
\end{table}

Algorithm~6 was implemented in software~\cite{Zhong}.

\section{Applications of Algorithm 6}\label{Shrikhande}

The Shrikhande graph and the Clebsch graph were used in \cite{KS} as a proof of concept for automatically generating proofs of non-word-representability for rather large graphs, for which finding proofs would otherwise be considered a challenging problem. A 10-line proof of non-word-representability of a certain minimal non-word-representable subgraph $C$ of the Clebsch graph using Algorithm 2 was produced in \cite{KS}. Algorithm 6 also produces a 10-line proof (both on $C$ and on the entire Clebsch graph); hence, we omit the discussion of the Clebsch graph here.

However, Algorithm 6 does shorten the 6-line proof given in \cite{KS} of non-word-representability of the minimal non-word-representable subgraph $S$ of the Shrikhande graph. Both of these graphs are presented in Figure~\ref{Shrikhande-pic}. By applying the same assumption on $S$ that vertex 1 is a source (as done in \cite{KS}), Algorithm~6 produces the following 4-line proof. \\

\begin{small}

\noindent
{\bf 1.} B2$\rightarrow$3 (Copy 2) O4$\rightarrow$3 (C1-2-3-4) O4$\rightarrow$5 (C1-3-4-5) O6$\rightarrow$5 (C1-4-5-6) B6$\rightarrow$9 (Copy 3) O5$\rightarrow$9 (C1-5-9-6) O8$\rightarrow$9 O4$\rightarrow$8 (C4-5-9-8) O3$\rightarrow$8 (C1-3-8-4) O7$\rightarrow$9 O2$\rightarrow$7 (C2-3-8-9-7) O7$\rightarrow$5 (C5-6-9-7) S:1-2-7-5

\noindent
{\bf 2.}  MC3 9$\rightarrow$6 O9$\rightarrow$7 O7$\rightarrow$5 (C5-6-9-7) O9$\rightarrow$5 (C1-5-9-6) O7$\rightarrow$2 (C1-2-7-5) O8$\rightarrow$3 O9$\rightarrow$8 (C2-3-8-9-7) O4$\rightarrow$8 (C4-5-9-8) S:1-4-8-3

\noindent
{\bf 3.} MC2 3$\rightarrow$2 O3$\rightarrow$4 (C1-2-3-4) O5$\rightarrow$4 (C1-3-4-5) O5$\rightarrow$6 (C1-4-5-6) B6$\rightarrow$9 (Copy 4) O7$\rightarrow$9 O5$\rightarrow$7 (C5-6-9-7) O5$\rightarrow$9 (C1-5-9-6) O2$\rightarrow$7 (C1-2-7-5) O3$\rightarrow$8 O8$\rightarrow$9 (C2-3-8-9-7) O8$\rightarrow$4 (C4-5-9-8) S:1-3-8-4

\noindent
{\bf 4.} MC4 9$\rightarrow$6 O9$\rightarrow$5 (C1-5-9-6) O9$\rightarrow$8 O8$\rightarrow$4 (C4-5-9-8) O8$\rightarrow$3 (C1-3-8-4) O9$\rightarrow$7 O7$\rightarrow$2 (C2-3-8-9-7) O5$\rightarrow$7 (C5-6-9-7) S:1-5-7-2. \\
\end{small}

As for the subgraph $S$ of the Kneser graph $K(8,3)$ in Figure~\ref{K83-fig}, for which no published proof of non-word-representability exists so far, one can use the software \cite{S} and assume that vertex 1 is a source to produce (in 25 seconds) a 136-line proof. However, with the same assumption, Algorithm~6 gives a much shorter proof of 110 lines, which we do not include in this paper because of its size (it would take about four pages); instead, we post the proof online at \cite{Zhong}.
 
\section{Concluding remarks}\label{concluding}

In this paper, we introduced an algorithm that significantly outperforms the best algorithm presented in \cite{KS} for automatically generating proofs of non-word-representability of graphs that can be verified by a human. The algorithm was implemented in publicly available software \cite{Zhong}, which substantially enlarges the set of graphs for which robust, human-verifiable proofs of non-word-representability can be produced. This advancement is expected to have a notable impact on future developments in the theory of word-representable graphs.

We leave it as an open problem to further improve Algorithm 6. A natural starting point for generating ideas is to analyze the graphs in Figure~\ref{15-non-word-repr} that require long proofs (for Algorithm 6, these are Graphs 1, 149, 818, and 5813) and to investigate whether different choices of branching vertices can reduce proof length. For example, for Graph 5813, is it possible to produce a proof shorter than 85 lines? This question is interesting in itself: constructing a shorter, ad-hoc proof for Graph 5813 could indicate that Algorithm 6 can indeed be improved.

We conclude by reiterating the following points from \cite{KS}: estimating the average efficiency of our algorithms, or of future algorithms, based on certain graph parameters (such as the number of cycles) is a challenging theoretical problem. Likewise, the time complexity of our algorithms presents an intriguing and challenging direction for research, which we have not addressed here because, as in \cite{KS}, our focus was on producing short proofs.


\end{document}